\documentclass[11pt]{article}
\usepackage[english]{babel}
\usepackage{indentfirst}
\usepackage{latexsym}
\usepackage[colorlinks=true, bookmarksnumbered=true, bookmarksopen=true,
bookmarksopenlevel=3, pdfstartview=FitH, linkcolor=cyan, pdfmenubar=true,
pdftoolbar=true, bookmarks=true,citecolor=cyan, urlcolor=magenta,
filecolor=cyan,menucolor=black,plainpages=false,pdfpagelabels]{hyperref}
\usepackage[lmargin=2cm,tmargin=2cm,bmargin=2cm,rmargin=2cm]{geometry}
\usepackage{amsbsy, amsmath,amsfonts,amssymb,amsthm}
\usepackage{times}
\usepackage{graphicx}
\usepackage{amsmath}
\usepackage{amsfonts}
\usepackage{amssymb}
\usepackage[english]{babel}
\usepackage[latin1]{inputenc}
\usepackage[T1]{fontenc}
\usepackage{amsmath,amssymb,graphics,amsthm}
\usepackage[usenames]{color}

\numberwithin{equation}{section}
\newtheorem{prop}{Proposition}[section]

\newtheorem{tma}[prop]{Theorem}

\newtheorem{obs}[prop]{Remark}
\newtheorem{lem}[prop]{Lemma}

\newtheorem{proposition}[prop]{Proposition}

\newcommand\reallywidehat[1]{\savestack{\tmpbox}{\stretchto{  \scaleto{    \scalerel*[\widthof{\ensuremath{#1}}]{\kern-.6pt\bigwedge\kern-.6pt}    {\rule[-\textheight/2]{1ex}{\textheight}}  }{\textheight}}{0.5ex}}\stackon[1pt]{#1}{\tmpbox}}
\def\fin { \vskip 0pt \hfill $\diamond$ \vskip 12pt}

\begin{document}

\title{On the 3D Euler equations with Coriolis force in borderline Besov
spaces}
\author{{Vladimir Angulo-Castillo$^{1}$}, \ {Lucas C. F. Ferreira$^{2}$}{%
\thanks{{Corresponding author. }\newline
{E-mail adresses: vladimirangulo01@gmail.com (V. Angulo-Castillo),
lcff@ime.unicamp.br (L.C.F. Ferreira).}\newline
{V. Angulo-Castillo was supported by CNPq, Brazil.\newline
LCF Ferreira was supported by FAPESP and CNPq, Brazil.}}} \\
{\small $^{1,2}$ Universidade Estadual de Campinas, Departamento de Matemá%
tica} \\
{\small CEP 13083-859, Campinas, SP, Brazil.}}
\date{}
\maketitle

\begin{abstract}
We consider the 3D Euler equations with Coriolis force (EC) in the whole
space. We show long-time solvability in Besov spaces for high speed of
rotation $\Omega $ and arbitrary initial data. For that, we obtain $\Omega $%
-uniform estimates and a blow-up criterion of BKM type in our framework. Our
initial data class is larger than previous ones considered for (EC) and
covers borderline cases of the regularity. The uniqueness of solutions is
also discussed. \medskip


{\small \medskip\noindent\textbf{Key:} Euler equations; Coriolis force;
Long-time solvability; Blow up; Besov-spaces}
\end{abstract}

\renewcommand{\abstractname}{Abstract}

\section{Introduction}

We consider the free incompressible Euler equations with Coriolis force
\begin{equation}
\left\{
\begin{split}
& \frac{\partial u}{\partial t}+\mathbb{P}\Omega e_{3}\times u+\mathbb{P}%
\left( u\cdot \nabla \right) u=0\text{ \ in }\ \mathbb{R}^{3}\times
(0,\infty ) \\
& \nabla \cdot u=0\text{ \ in }\mathbb{R}^{3}\times (0,\infty ) \\
& u(x,0)=u_{0}(x)\text{ \ in }\ \mathbb{R}^{3}
\end{split}%
\right. ,  \label{prin-problem-leray}
\end{equation}%
where $u(x,t)=\left( u_{1}(x,t),u_{2}(x,t),u_{3}(x,t)\right) $ stands for
the velocity field, $\mathbb{P}=(\delta _{jk}+R_{j}R_{k})_{1\leq j,k\leq 3}$
is the Leray-Helmholtz projection and $R_{j}$ denotes the $j$-th Riesz
transform. The Coriolis parameter $\Omega \in \mathbb{R}$ corresponds to
twice the speed of rotation around the vertical unit vector $e_{3}=(0,0,1).$
The initial velocity is denoted by $u_{0}=u_{0}(x)=\left(
u_{0,1}(x),u_{0,2}(x),u_{0,3}(x)\right) $ and satisfies the compatibility
condition $\nabla \cdot u_{0}=0$. The reader is referred to the book \cite%
{Chemin-Gallagher} for more details about the physical model. Throughout the
paper, we denote spaces of scalar and vector functions abusively in the same
way; for example, we write $u_{0}\in H^{s}(\mathbb{R}^{3})$ instead of $%
u_{0}\in (H^{s}(\mathbb{R}^{3}))^{3}.$

The system (\ref{prin-problem-leray}) has been studied by several authors in
the case $\Omega =0$ that corresponds to the classical Euler equations (E).
In what follows we give a brief review of some of these results. In the
framework of Sobolev spaces, Kato \cite{Kato1972} showed that (E) has a
unique local-in-time solution $u\in C\left( \left[ 0,T\right] ;H^{s}\left(
\mathbb{R}^{3}\right) \right) \cap C^{1}(\left[ 0,T\right] ;H^{s-1}\left(
\mathbb{R}^{3}\right) )$ for $u_{0}\in H^{s}\left( \mathbb{R}^{3}\right) $
with an integer $s\geq 3$ where $T=T(\left\Vert u_{0}\right\Vert _{H^{s}(%
\mathbb{R}^{3})})$. In \cite{Kato-ibero-1986}, Kato and Ponce proved that if
$s>2/p+1$, $1<p<\infty $ and $u_{0}\in H_{p}^{s}\left( \mathbb{R}^{2}\right)
,$ then there exists a unique 2D global solution $u\in C\left( \left[
0,\infty \right) ;H_{p}^{s}\left( \mathbb{R}^{2}\right) \right) $. Later, in
\cite{KatoPonce1988} they considered $n\geq 2$ and proved that for $s>n/p+1,$
$1<p<\infty $ and $u_{0}\in H_{p}^{s}\left( \mathbb{R}^{n}\right) ,$ there
exist $T>0$ and a unique solution $u\in C\left( \left[ 0,T\right]
;H_{p}^{s}\left( \mathbb{R}^{n}\right) \right) \cap C^{1}\left( \left[ 0,T%
\right] ;H_{p}^{s-1}\left( \mathbb{R}^{n}\right) \right) .$ Temam \cite%
{Temam1975} extended the results of Kato \cite{Kato1972} to $H^{m}$ and $%
W^{m,p}$ in bounded domains (see also Ebin-Marsden \cite{Marsden-1} and
Bourguignon-Brezis \cite{Brezis}). For existence and uniqueness results in
Holder $C^{k,\gamma }$ and Triebel-Lizorkin $F_{p,q}^{s}$ spaces, the reader
is referred to \cite{Chemin-1} and \cite{Chae-2, Chae-3}, respectively.

In the context of Besov spaces, Chae \cite{Chae2004} and Zhou \cite{Zhou-2}
proved that (E) has a unique solution $u\in C([0,T];\linebreak{%
B_{p,1}^{n/p+1}\left( \mathbb{R}^{n}\right)})$ for $1<p<\infty $ and $n\geq
3 $ (see also \cite{VS} for $n=2$). After, the borderline cases $p=\infty $
\cite{PP1} and $p=1$ \cite{PP2} was considered by Pak and Park. Takada \cite%
{Takada2008} showed existence-uniqueness in Besov type spaces based on weak-$%
L^{p}$ with $1<p<\infty $ and $n\geq 3.$ The exponent $s=\frac{n}{p}+1$ is
critical for (E) in $H_{p}^{s}$ and $B_{p,q}^{s}$-spaces. In fact, Bourgain
and Li \cite{Bourgain} showed that (E) is ill-posed in $H_{p}^{s}$ and $%
B_{p,q}^{n/p+1}$ for $1\leq p<\infty $, $1<q\leq \infty $ and $n=2,3$. So,
it is natural to consider $q=1$ when $s=\frac{n}{p}+1$. The critical case is
also of special interest because the regularity index $s-1=\frac{n}{p}$ of
the vorticity $\nabla \times u$ corresponds to a critical case of Sobolev
type embeddings. Motivated by the symbol $>,$ the case $s>\frac{n}{p}+1$ has
been named in the literature as supercritical.

For $\Omega \neq 0$, Dutrifoy \cite{Dutrifoy2005} showed long-time existence
of solutions for (\ref{prin-problem-leray}) with lower bound on the
existence-time $T_{\Omega }\gtrsim \log {\log {|\Omega |}}$ provided that $%
|\Omega |$ is large enough and $u_{0}$ belongs to a certain Sobolev type
class. Also, Dutrifoy \cite{Dutrifoy2004} and Charve \cite{Charve2006}
obtained analogous results for quasigeostrophic systems. Recently, for $%
s>s_{0}=\frac{3}{2}+1$ and $u_{0}\in H^{s}(\mathbb{R}^{3}),$ Koh, Lee and
Takada \cite{Koh2015} proved that there exists a unique local in time
solution $u$ for (\ref{prin-problem-leray}) in the class $C([0,T];H^{s}(%
\mathbb{R}^{3}))\cap C^{1}([0,T];H^{s-1}(\mathbb{R}^{3})).$ Moreover,
assuming that $s>s_{1}=\frac{5}{2}+1,$ they showed that their solutions can
be extended to long-time intervals $[0,T_{\Omega }]$ provided that the speed
of rotation is large enough. For the viscous case, we refer the reader to
the works \cite{Babin1,Babin2,Chemin-Gallagher,Takada} for global
well-posedness in Sobolev spaces with $\left\vert \Omega \right\vert $ large
enough and to the papers \cite{Giga-1,Takada2} (and their references) for
results about global well-posedness with $\Omega $-uniform smallness
condition on initial data in different types of critical spaces (e.g., in
Fourier Besov spaces).

In view of the previous results for (\ref{prin-problem-leray}) and (E), it
is natural to wonder about the borderline cases $s_{0}$ and $s_{1}.$ In this
paper we extend the results of \cite{Koh2015} by treating these two cases in
the framework of Besov spaces. To be more precise, we consider the critical
regularity $s_{0}$ and show local-in-time existence and uniqueness of
solutions for initial data in the critical Besov space $B_{2,1}^{s_{0}}$
with smallness condition on the existence-time uniformly in $\Omega \in
\mathbb{R}.$ After, for large Coriolis parameter $\left\vert \Omega
\right\vert ,$ we obtain long-time solvability of (\ref{prin-problem-leray})
in $B_{2,1}^{s}$ in the borderline case $s=s_{1}$. \ It is worth to observe
that $H^{s}\subset B_{2,1}^{5/2}$ and \ $H^{s}\subset B_{2,1}^{7/2}$ for $s>%
\frac{3}{2}+1$ and $s>\frac{5}{2}+1,$ respectively, and so our result
provides a larger class for both local and long time solvability of (\ref%
{prin-problem-leray}).

Our main result reads as follows.

\begin{tma}
\label{principal-theorem}

\begin{enumerate}
\item[$(i)$] Let $u_{0}\in B_{2,1}^{5/2}(\mathbb{R}^{3})$ satisfy $\nabla
\cdot u_{0}=0$. There exists $T=T(\Vert u_{0}\Vert _{B_{2,1}^{5/2}})>0$ such
that (\ref{principal-theorem}) has a unique solution $u\in
C([0,T];B_{2,1}^{5/2}(\mathbb{R}^{3}))\cap C^{1}([0,T];B_{2,1}^{3/2}(\mathbb{%
R}^{3}))$, for all $\Omega \in \mathbb{R}$.

\item[$(ii)$] Let $0<T<\infty $ and $u_{0}\in B_{2,1}^{7/2}(\mathbb{R}^{3})$
be such that $\nabla \cdot u_{0}=0$. There exists $\Omega _{0}=\Omega
_{0}(T,\Vert u_{0}\Vert _{B_{2,1}^{7/2}})>0$ such that (\ref%
{principal-theorem}) has a unique solution $u\in C([0,T];B_{2,1}^{7/2}(%
\mathbb{R}^{3}))\cap C^{1}([0,T];B_{2,1}^{5/2}(\mathbb{R}^{3}))$ provided
that $|\Omega |\geq \Omega _{0}.$
\end{enumerate}
\end{tma}

Considering $\Omega =0,$ item $(i)$ recovers the local existence result by
Chae \cite{Chae2004} and Zhou \cite{Zhou-2} for Euler equations in $%
B_{p,1}^{n/p+1}(\mathbb{R}^{n})$ in the case $p=2$ and $n=3$. Assuming
further regularity on the initial data, item $(ii)$ shows that local
solutions can be extended to arbitrary large time $T>0$ provided that $%
\left\vert \Omega \right\vert $ is large enough and so it resembles results
for the 2D Euler equations (see \cite{VS, Chae2004}). In fact, we recall
that existence of smooth solutions for the 3D Euler equations is an
outstanding open problem. Long-time solvability type results for (\ref%
{prin-problem-leray}) with arbitrary data show a smoothing effect connected
to the speed of rotation $\Omega $ (see \cite{Chemin-Gallagher}).

Finally, we comment on some technical points in our results. The general
strategy of this paper consists in three basic steps: approximation scheme;
\textit{a priori} $\Omega $-uniform estimates and passing to the limit for obtaining
local-in time solutions; blow-up criterion and long-time solvability. This
is the same one employed by \cite{Koh2015} in $H^{s}$-spaces however here we
need to carry out the necessary estimates in the borderline Besov spaces $%
B_{2,1}^{s_{0}}$ and $B_{2,1}^{s_{1}}$. In order to pass the limit in the
approximation scheme $\{u^{\delta }\}_{\delta >0}$, the authors of \cite%
{Koh2015} relied on the Hilbert structure of $H^{s}$-spaces. Since our
setting has not such property, we need to control $u^{\delta }$ by means of estimates involving localization and $B_{2,1}^{s}$-norms (see, e.g., Lemma \ref{est-seq-apro-0},
Proposition \ref{est-seq-apro-1} and proof of Theorem \ref{principal-theorem}%
). In order to cover the endpoints $s_{0}$ and $s_{1}$ of the ranges in \cite%
{Koh2015}, we are inspired by previous results for the Euler equations (E)
\cite{Chae2004, PP1, Zhou-2} and consider $B_{2,1}^{s}$-spaces (and the
embedding $B_{2,1}^{3/2}\hookrightarrow L^{\infty }$ in $\mathbb{R}^{3}$)
that allow us to have $\nabla u\in L^{\infty }$ for $s=s_{0}$ (which is not true in $H^{s}$) and control globally in time $U(t)=\int_{0}^{t}\Vert \nabla u(\tau)\Vert _{L^{\infty }}\ d\tau$ for large $\Omega$ when $s=s_{1}$. The quantity $U(t)$ is used to derive a blow-up criterion
and obtain long-time solutions. Also, we show Lemma \ref{continuous} that
deals with the time-continuity of weak solutions for (\ref%
{prin-problem-leray}) and is useful to prove time-regularity of solutions
obtained as limit of the approximation scheme. For that matter, we extend
\cite[Lemma 2.1]{PP2} (that considered solutions of (E) in $B_{1,1}^{d+1}(%
\mathbb{R}^{d})$) to the Euler Coriolis equations in $B_{p,1}^{3/p+1}(%
\mathbb{R}^{3})$ with $1\leq p<\infty .$

The plan of this paper is as follows. The next section is devoted to some
preliminaries about product and commutator estimates in Besov spaces and
projection operators linked to the Coriolis term. In Section 3, we deal with
the approximation scheme $\{u^{\delta }\}_{\delta >0}$ and show local
existence on $[0,T]$ with $T>0$ independent of $\delta $ and $\Omega .$ The
proof of Theorem \ref{principal-theorem} is given in Section 4 through three
subsections: item $(i)$ in subsection 4.1, blow-up criterion in subsection
4.2, and item $(ii)$ in subsection 4.3.

\section{Function spaces and projection operators}

This section is devoted to some preliminaries about Besov spaces. We refer
the reader to \cite{Berg-Lofstrom} for more details on these spaces and
their properties. Also, we recall two projection operators that will be
useful for our purposes.

Let $\mathcal{S}(\mathbb{R}^{3})$ and $\mathcal{S}^{\prime }(\mathbb{R}^{3})$
stand for the Schwartz class and the space of tempered distributions,
respectively. Let $\widehat{f}$ denote the Fourier transform of $f\in
\mathcal{S}^{\prime }$. Consider a nonnegative radial function $\phi _{0}\in
\mathcal{S}(\mathbb{R}^{3})$ satisfying $0\leq \widehat{\phi }_{0}(\xi )\leq
1$ for all $\xi \in \mathbb{R}^{3}$, $\mbox{supp}\ \widehat{\phi }%
_{0}\subset \{\xi \in \mathbb{R}^{3}:\frac{1}{2}\leq |\xi |\leq 2\}$ and
\begin{equation*}
\sum_{j\in \mathbb{Z}}\widehat{\phi }_{j}(\xi )=1\ \ \mbox{for all}\ \ \xi
\in \mathbb{R}^{3}\backslash \{0\},
\end{equation*}%
where $\phi _{j}(x):=2^{3j}\phi _{0}(2^{j}x).$ For $k\in \mathbb{Z}$, we
define the function $S_{k}\in \mathcal{S}$ as
\begin{equation*}
\widehat{S}_{k}(\xi )=1-\sum_{j\geq k+1}\widehat{\phi }_{j}(\xi )
\end{equation*}%
and denote $\psi =S_{0}$. For $f\in \mathcal{S}^{\prime }(\mathbb{R}^{3}),$
the Littlewood-Paley operator $\Delta _{j}$ is defined by $\Delta
_{j}f:=\phi _{j}\ast f.$

Let $s\in \mathbb{R}$ and $1\leq p,q\leq \infty $ and let $\mathcal{P}$
denote the set of polynomials with $3$ variables. The homogeneous Besov
space ${\dot{B}}_{p,q}^{s}(\mathbb{R}^{3})$ is the set of all $f\in \mathcal{%
S}^{\prime }(\mathbb{R}^{3})/\mathcal{P}$ such that
\begin{equation*}
\Vert f\Vert _{{\dot{B}}_{p,q}^{s}}:=\Vert \{2^{sj}\Vert \Delta _{j}f\Vert
_{L^{p}}\}_{j\in \mathbb{Z}}\Vert _{l^{q}(\mathbb{Z})}<\infty .
\end{equation*}%
The inhomogeneous version of ${\dot{B}}_{p,q}^{s}$, denoted by $B_{p,q}^{s}(%
\mathbb{R}^{3})$, is defined as the set of all $f\in \mathcal{S}^{\prime }(%
\mathbb{R}^{3})$ such that
\begin{equation*}
\Vert f\Vert _{B_{p,q}^{s}}:=\Vert \{2^{sj}\Vert \Delta _{j}f\Vert
_{L^{p}}\}_{j\in \mathbb{N}}\Vert _{l^{q}(\mathbb{N})}+\Vert \psi \ast
f\Vert _{L^{p}}.
\end{equation*}%
The pairs $({\dot{B}}_{p,q}^{s},\Vert \cdot \Vert _{{\dot{B}}_{p,q}^{s}})$
and $(B_{p,q}^{s},\Vert \cdot \Vert _{B_{p,q}^{s}})$ are Banach spaces. For $%
s>0$, we have the equivalence
\begin{equation}
\Vert f\Vert _{B_{p,q}^{s}}\sim \Vert f\Vert _{{\dot{B}}_{p,q}^{s}}+\Vert
f\Vert _{L^{p}}.  \label{normequivalence}
\end{equation}

\begin{lem}[Bernstein inequality]
\label{inequalitybernstein} Assume that $f\in L^{p}$, $1\leq p\leq \infty $,
and $\mbox{supp}\ \widehat{f}\subset \{\xi \in \mathbb{R}^{3}:2^{j-2}\leq
|\xi |<2^{j}\}$. Then there exists a constant $C=C(k)>0$ such that
\begin{equation*}
C^{-1}2^{jk}\Vert f\Vert _{L^{p}}\leq \Vert D^{k}f\Vert _{L^{p}}\leq
C2^{jk}\Vert f\Vert _{L^{p}}.
\end{equation*}
\end{lem}

\begin{obs}
\label{remark1}As a consequence of the above lemma we have the following
equivalence
\begin{equation}
\Vert D^{k}f\Vert _{{\dot{B}}_{p,q}^{s}}\sim \Vert f\Vert _{{\dot{B}}%
_{p,q}^{s+k}}.  \label{equiv-1}
\end{equation}
We also recall the estimate (see, e.g., \cite{Takada2008})
\begin{equation}
\Vert f\Vert _{L^{\infty }}\leq C\Vert f\Vert _{B_{p,q}^{s}},
\label{est-aux-1}
\end{equation}%
where $s>n/p$ with $1\leq p,q\leq \infty $, or $s=n/p$ with $1\leq p\leq
\infty $ and $q=1$. Thus, for $s>n/p+1$ with $1\leq p,q\leq \infty $ or $%
s=n/p+1$ with $1\leq p\leq \infty $ and $q=1$, we have the estimates
\begin{equation}
\Vert \nabla f\Vert _{L^{\infty }}\leq \Vert \nabla f\Vert
_{B_{p,q}^{s-1}}\leq \Vert f\Vert _{B_{p,q}^{s}}.  \label{est-aux-2}
\end{equation}
\end{obs}

The following lemma contains product estimates in the framework of Besov
spaces (see \cite{Chae2004}).

\begin{lem}
\label{inequalityholder} Let $s>0$, $1\leq p,q\leq \infty $, $1\leq
p_{1},p_{2}\leq \infty $ and $1\leq r_{1},r_{2}\leq \infty $ satisfy $\frac{1%
}{p}=\frac{1}{p_{1}}+\frac{1}{p_{2}}=\frac{1}{r_{1}}+\frac{1}{r_{2}}.$ Then
there exists a universal constant $C>0$ such that%
\begin{eqnarray*}
\Vert fg\Vert _{{\dot{B}}_{p,q}^{s}} &\leq &C(\Vert f\Vert _{{\dot{B}}%
_{p_{1},q}^{s}}\Vert g\Vert _{L^{p_{2}}}+\Vert g\Vert _{{\dot{B}}%
_{r_{1},q}^{s}}\Vert f\Vert _{L^{r_{2}}}) \\
\Vert fg\Vert _{B_{p,q}^{s}} &\leq &C(\Vert f\Vert _{B_{p_{1},q}^{s}}\Vert
g\Vert _{L^{p_{2}}}+\Vert g\Vert _{B_{r_{1},q}^{s}}\Vert f\Vert
_{L^{r_{2}}}).
\end{eqnarray*}
\end{lem}

In the next two lemmas we recall estimates in ${\dot{B}}_{p,q}^{s}$ and $%
B_{p,q}^{s}$ for the commutator (see \cite{Chae2004, Takada2008})
\begin{equation*}
\lbrack v\cdot \nabla ,\Delta _{j}]u=v\cdot \nabla (\Delta _{j}u)-\Delta
_{j}(v\cdot \nabla u).
\end{equation*}

\begin{lem}
\label{2.4} Let $1<p<\infty $ and $1\leq q\leq \infty $.

\begin{enumerate}
\item[$(i)$] Let $s>0$, $v\in {\dot{B}}_{p,q}^{s}(\mathbb{R}^{n})$ with $%
\nabla v\in L^{\infty }(\mathbb{R}^{n})$ and $\nabla \cdot v=0,$ and $\theta
\in {\dot{B}}_{p,q}^{s}(\mathbb{R}^{n})$ with $\nabla \theta \in L^{\infty }(%
\mathbb{R}^{n})$. Then, there exists a universal constant $C>0$ such that
\begin{equation*}
\left( \sum_{j\in \mathbb{Z}}2^{sjq}\Vert \lbrack v\cdot \nabla ,\Delta
_{j}]\theta \Vert _{L^{p}}^{q}\right) ^{1/q}\leq C\left( \Vert \nabla v\Vert
_{L^{\infty }}\Vert \theta \Vert _{{\ \dot{B}}_{p,q}^{s}}+\Vert \nabla
\theta \Vert _{L^{\infty }}\Vert v\Vert _{{\dot{B}}_{p,q}^{s}}\right) .
\end{equation*}

\item[$(ii)$] Let $s>-1$, $v\in {\dot{B}}_{p,q}^{s+1}(\mathbb{R}^{n})$ with $%
\nabla v\in L^{\infty }(\mathbb{R}^{n})$ and $\nabla \cdot v=0,$ and $\theta
\in {\dot{B}}_{p,q}^{s}(\mathbb{R}^{n})\cap L^{\infty }(\mathbb{R}^{n})$.
Then, there exists a universal constant $C>0$ such that
\begin{equation*}
\left( \sum_{j\in \mathbb{Z}}2^{sjq}\Vert \lbrack v\cdot \nabla ,\Delta
_{j}]\theta \Vert _{L^{p}}^{q}\right) ^{1/q}\leq C\left( \Vert \nabla v\Vert
_{L^{\infty }}\Vert \theta \Vert _{{\ \dot{B}}_{p,q}^{s}}+\Vert \theta \Vert
_{L^{\infty }}\Vert v\Vert _{{\dot{B}}_{p,q}^{s+1}}\right) .
\end{equation*}
\end{enumerate}
\end{lem}

\begin{lem}
\label{bilinearest2} Let $1<p<\infty $ and let $s>3/p+1$ with $1\leq q\leq
\infty $ or $s=3/p+1$ with $q=1$. Then, there exists a constant $C>0$ such
that
\begin{equation*}
\left( \sum_{j\in \mathbb{Z}}2^{jqs}\Vert (S_{j-2}u\cdot \nabla )\Delta
_{j}u-\Delta _{j}(u\cdot \nabla )u\Vert _{L^{p}}^{q}\right) ^{1/q}\leq
C\Vert \nabla u\Vert _{L^{\infty }}\Vert u\Vert _{B_{p,q}^{s}},
\end{equation*}%
for all $u\in B_{p,q}^{s}(\mathbb{R}^{3})$ with $\nabla \cdot u=0$.
\end{lem}

\bigskip In order to handle the Coriolis term, we will need the following
projection operators $P_{\pm }:L^{2}(\mathbb{R}^{3})^{3}\longrightarrow
L^{2}(\mathbb{R}^{3})^{3}$ given by%
\begin{equation*}
P_{\pm }v:=\frac{1}{2}\left( \mathbb{P}v\pm i\frac{D}{|D|}\times v\right) ,
\end{equation*}%
where $\frac{D}{|D|}\times $ is defined by means of the Fourier transform as
$(\frac{D}{|D|}\times v)^{\widehat{}}(\xi ):=\frac{\xi }{|\xi |}\times
\widehat{v}(\xi )$.

The next lemma contains basic properties of $P_{\pm }$ and can be found in
\cite{Dutrifoy2005, Koh2015}.

\begin{lem}
\label{projection}The projections $P_{\pm }$ satisfy $P_{\pm }\mathbb{P}%
=P_{\pm }$. Moreover, if $\ \nabla \cdot v=0$ we have that $v=P_{+}v+P_{-}v$%
, $\mathbb{P}\left( e_{3}\times v\right) =-i\frac{D_{3}}{|D|}(P_{+}v-P_{-}v)$%
, $P_{\pm }P_{\pm }=P_{\pm }$, and $P_{\pm }P_{\mp }=0$.
\end{lem}

\section{Approximation scheme}

Let $u_{0}$ be the initial velocity in (\ref{prin-problem-leray}). For $%
0<\delta <1$, we consider the approximate parabolic problem
\begin{equation}
\left\{
\begin{split}
& \frac{\partial u^{\delta }}{\partial t}-\delta \Delta u^{\delta }+\mathbb{P%
}\Omega e_{3}\times u^{\delta }+\mathbb{P}\left( u^{\delta }\cdot \nabla
\right) u^{\delta }=0\text{ in }\mathbb{R}^{3}\times (0,\infty ), \\
& \ \nabla \cdot u^{\delta }=0\text{ in }\mathbb{R}^{3}\times (0,\infty ), \\
& u^{\delta }(x,0)=u_{0}(x)\text{ in }\mathbb{R}^{3}.
\end{split}%
\right.  \label{Approx-scheme}
\end{equation}%
We are going to show that the above problem has a solution for each $\delta
>0$ in a suitable class involving Besov spaces. For that matter, first we
recall some estimates for the heat semigroup $\{e^{t\Delta }\}_{t\geq 0}$ in
$B_{p,q}^{s}$ (see, e.g., \cite{Kozono2003}).

\begin{lem}
\label{kozono} Let $s_{0}\leq s_{1}$ and $1\leq p,q\leq \infty $. Then there
exists a constant $C>0$ (independent of $p,q$ and $t>0$) such that
\begin{equation*}
\Vert e^{t\Delta }f\Vert _{B_{p,q}^{s_{1}}}\leq C(1+t^{-\frac{1}{2}%
(s_{1}-s_{0})})\Vert f\Vert _{B_{p,q}^{s_{0}}},
\end{equation*}%
for all $f\in {B_{p,q}^{s_{0}}}(\mathbb{R}^{3})$.
\end{lem}

We start by showing estimates for the bilinear term of the mild formulation
for (\ref{Approx-scheme}).

\begin{lem}
\label{3.5} Let $0<\delta <1$ and $1<p<\infty $.

\begin{enumerate}
\item[$(i)$] There exists $C>0$ such that
\begin{equation}
\sup_{0<t<T}\left\Vert \int_{0}^{t}e^{\delta (t-\tau )\Delta }\mathbb{P}%
\left( u(\tau )\cdot \nabla \right) v(\tau )\ d\tau \right\Vert _{{%
B_{p,1}^{3/p}}}\leq C\sup_{0<t<T}\Vert u(t)\Vert _{B_{p,1}^{3/p}}\Vert
v\Vert _{L^{1}(0,T;{B_{p,1}^{3/p+1}})},  \label{aux-est-heat-10}
\end{equation}%
for all $u\in C([0,T];{B_{p,1}^{3/p}}(\mathbb{R}^{3}))$ and $v\in L^{1}(0,T;{%
B_{p,1}^{3/p+1}}(\mathbb{R}^{3})).$

\item[$(ii)$] Let $k=1,2.$ There exists $C>0$ such that
\begin{equation*}
\left\Vert \int_{0}^{t}e^{\delta (t-\tau )\Delta }\mathbb{P}\left( u(\tau
)\cdot \nabla \right) v(\tau )\ d\tau \right\Vert _{L^{1}(0,T;{%
B_{p,1}^{3/p+k}})}\leq C\left( T+T^{\frac{1}{2}}\delta ^{-\frac{1}{2}%
}\right) \sup_{0<t<T}\Vert u(t)\Vert _{B_{p,1}^{3/p+k-1}}\Vert v\Vert
_{L^{1}(0,T;{B_{p,1}^{3/p+k}})},
\end{equation*}%
for all $u\in C([0,T];{B_{p,1}^{3/p}}(\mathbb{R}^{3}))$ with $\nabla \cdot
u=0$ and $v\in L^{1}(0,T;{B_{p,1}^{3/p+k}}(\mathbb{R}^{3})).$
\end{enumerate}
\end{lem}

\textbf{Proof.} For $1<p<\infty ,$ we have that $\Vert e^{\delta t\Delta
}\Vert _{B_{p,1}^{3/p}\rightarrow B_{p,1}^{3/p}}\leq 1$ and $\mathbb{P}$ is
bounded in $B_{p,1}^{3/p}.$ So, we can estimate
\begin{equation*}
\left\Vert \int_{0}^{t}e^{\delta (t-\tau )\Delta }\mathbb{P}\left( u(\tau
)\cdot \nabla \right) v(\tau )\ d\tau \right\Vert _{{B_{p,1}^{3/p}}}\leq
C\int_{0}^{t}\Vert \left( u(\tau )\cdot \nabla \right) v(\tau )\Vert
_{B_{p,1}^{3/p}}\ d\tau .
\end{equation*}%
From Lemmas \ref{inequalityholder} and \ref{inequalitybernstein}, it follows
that
\begin{equation*}
\Vert \left( u(\tau )\cdot \nabla \right) v(\tau )\Vert _{B_{p,1}^{3/p}}\leq
C\Vert u(\tau )\Vert _{B_{p,1}^{3/p}}\Vert v(\tau )\Vert _{B_{p,1}^{3/p+1}}
\end{equation*}%
and then
\begin{equation*}
\begin{split}
\left\Vert \int_{0}^{t}e^{\delta (t-\tau )\Delta }\mathbb{P}\left( u(\tau
)\cdot \nabla \right) v(\tau )\ d\tau \right\Vert _{{B_{p,1}^{3/p}}}& \leq
C\int_{0}^{t}\Vert u(\tau )\Vert _{B_{p,1}^{3/p}}\Vert v(\tau )\Vert
_{B_{p,1}^{3/p+1}}\ d\tau \\
& \leq C\sup_{0<t<T}\Vert u(t)\Vert _{B_{p,1}^{3/p}}\Vert v\Vert _{L^{1}(0,T;%
{B_{p,1}^{3/p+1}})},
\end{split}%
\end{equation*}%
for all $0<t<T$, which gives (\ref{aux-est-heat-10}).

By Minkowski inequality and Lemmas \ref{kozono}, \ref{inequalityholder}
and \ref{inequalitybernstein}, we have that

\begin{equation}
\begin{split}
\Vert \int_{0}^{t}& e^{\delta (t-\tau )\Delta }\mathbb{P}\left( u(\tau
)\cdot \nabla \right) v(\tau )\ d\tau \Vert _{{B_{p,1}^{3/p+k}}} \\
& \leq C\int_{0}^{t}\{1+\delta ^{-\frac{1}{2}}(t-\tau )^{-\frac{1}{2}%
}\}\left( \Vert u(\tau )\Vert _{B_{p,1}^{3/p}}\Vert v(\tau )\Vert
_{B_{p,1}^{3/p+k}}+\Vert u(\tau )\Vert _{B_{p,1}^{3/p+k-1}}\Vert v(\tau
)\Vert _{B_{p,1}^{3/p+1}}\right) \ d\tau \\
& \leq C\int_{0}^{t}\{1+\delta ^{-\frac{1}{2}}(t-\tau )^{-\frac{1}{2}%
}\}\Vert u(\tau )\Vert _{B_{p,1}^{3/p+k-1}}\Vert v(\tau )\Vert
_{B_{p,1}^{3/p+k}}\ d\tau \\
& \leq C\sup_{0<t<T}\Vert u(t)\Vert _{B_{p,1}^{3/p+k-1}}\left\{ \Vert v\Vert
_{L^{1}(0,T;{B_{p,1}^{3/p+k}})}+\delta ^{-\frac{1}{2}}\int_{0}^{t}(t-\tau
)^{-\frac{1}{2}}\Vert v(\tau )\Vert _{B_{p,1}^{3/p+k}}\ d\tau \right\} ,
\end{split}
\label{aux-est-mild-1}
\end{equation}
for all $0<t<T$ and $k=1,2$. We can now compute the norm $\left\Vert \cdot
\right\Vert _{L^{1}(0,T)}$ in (\ref{aux-est-mild-1}) to obtain
\begin{equation*}
\begin{split}
\Vert \int_{0}^{t}& e^{\delta (t-\tau )\Delta }\mathbb{P}\left( u(\tau
)\cdot \nabla \right) v(\tau )\ d\tau \Vert _{L^{1}(0,T;{B_{p,1}^{3/p+k}})}
\\
& \leq C\sup_{0<t<T}\Vert u(t)\Vert _{B_{p,1}^{3/p+k-1}}\left\{ T\Vert
v\Vert _{L^{1}(0,T;{B_{p,1}^{3/p+k}})}+\delta ^{-\frac{1}{2}%
}\int_{0}^{T}\Vert v(\tau )\Vert _{B_{p,1}^{3/p+k}}\int_{\tau }^{T}(t-\tau
)^{-\frac{1}{2}}dt\ d\tau \right\} \\
& \leq C\left( T+T^{\frac{1}{2}}\delta ^{-\frac{1}{2}}\right)
\sup_{0<t<T}\Vert u(t)\Vert _{B_{p,1}^{3/p+k-1}}\Vert v\Vert _{L^{1}(0,T;{%
B_{p,1}^{3/p+k}})}.
\end{split}%
\end{equation*}%
\fin

Before proceeding, we recall that $AC([0,T];X)$ denotes the set of all $X$%
-valued absolutely continuous functions on $[0,T]$. The next lemma ensures
the existence of strong solution for (\ref{Approx-scheme}). The proof
follows essentially the same steps of \cite[Lemma 3.1.]{Koh2015} but using
estimates in Besov spaces instead of Sobolev spaces.

\begin{lem}
\label{est-seq-apro-0}Let $1<p<\infty $, $\delta \in (0,1)$ and $\Omega \in
\mathbb{R}$. Assume that $u_{0}\in B_{p,1}^{3/p+1}(\mathbb{R}^{3})$ and $%
\nabla \cdot u_{0}=0$. Then there exists a positive time $T_{\delta ,\Omega
}=T(\delta ,|\Omega |,\Vert u_{0}\Vert _{B_{p,1}^{3/p+1}})$ such that (\ref%
{Approx-scheme}) has a unique strong solution $u^{\delta }$ satisfying
\begin{equation}
u^{\delta }\in C([0,T_{\delta ,\Omega }];{B_{p,1}^{3/p+1}}(\mathbb{R}%
^{3}))\cap AC([0,T_{\delta ,\Omega }];{B_{p,1}^{3/p}}(\mathbb{R}^{3}))\cap
L^{1}(0,T_{\delta ,\Omega };{B_{p,1}^{3/p+2}}(\mathbb{R}^{3}))  \label{eq3.1}
\end{equation}
\end{lem}

\textbf{Proof.} \ Firstly, we consider the mild formulation for (\ref%
{Approx-scheme})
\begin{equation}
u^{\delta }(t)=e^{\delta t\Delta }u_{0}-\int_{0}^{t}e^{\delta (t-\tau
)\Delta }\mathbb{P}\Omega e_{3}\times u^{\delta }(\tau )\ d\tau
-\int_{0}^{t}e^{\delta (t-\tau )\Delta }\mathbb{P}\left( u^{\delta }(\tau
)\cdot \nabla \right) u^{\delta }(\tau )\ d\tau  \label{eq3.2}
\end{equation}%
and show the existence of a local in time solution. Lemma \ref{kozono}
yields the estimate
\begin{equation*}
\Vert e^{\delta t\Delta }f\Vert _{L^{1}(0,T;{B_{p,1}^{3/p+2}})}\leq C(T+T^{%
\frac{1}{2}}\delta ^{-\frac{1}{2}})\Vert f\Vert _{B_{p,1}^{3/p+1}},
\end{equation*}%
for all $f\in {B_{p,1}^{3/p+1}}$. Thus, for all $0<T<\infty $ we have
\begin{equation}
\sup_{0\leq t\leq T}\Vert e^{\delta t\Delta }u_{0}\Vert
_{B_{p,1}^{3/p+1}}+L_{\delta ,T}^{-1}\Vert e^{\delta t\Delta }u_{0}\Vert
_{L^{1}(0,T;B_{p,1}^{3/p+2})}\leq C_{0}\Vert u_{0}\Vert _{B_{p,1}^{3/p+1}},
\label{eq3.3}
\end{equation}%
where $C_{0}>0$ is a constant and $L_{\delta ,T}=(T+T^{\frac{1}{2}}\delta ^{-%
\frac{1}{2}}).$

Consider the map
\begin{equation*}
\mathcal{B}(u^{\delta })(t)=e^{\delta t\Delta }u_{0}-\int_{0}^{t}e^{\delta
(t-\tau )\Delta }\mathbb{P}\Omega e_{3}\times u^{\delta }(\tau )\ d\tau
-\int_{0}^{t}e^{\delta (t-\tau )\Delta }\mathbb{P}\left( u^{\delta }(\tau
)\cdot \nabla \right) u^{\delta }(\tau )\ d\tau
\end{equation*}%
and the complete metric space
\begin{equation*}
Z_{T}:=\left\{ u\in C([0,T];B_{p,1}^{3/p+1}(\mathbb{R}^{3}))\cap
L^{1}(0,T;B_{p,1}^{3/p+2}(\mathbb{R}^{3}));\text{ }\nabla \cdot u=0\text{
and }\Vert u\Vert _{Z_{T}}\leq 2C_{0}\Vert u_{0}\Vert
_{B_{p,1}^{3/p+1}}\right\}
\end{equation*}%
whose norm is given by
\begin{equation*}
\Vert u\Vert _{Z_{T}}:=\sup_{0\leq t\leq T}\Vert u(t)\Vert
_{B_{p,1}^{3/p+1}}+L_{\delta ,T}^{-1}\Vert u\Vert
_{L^{1}(0,T;B_{p,1}^{3/p+2})}.
\end{equation*}%
We claim that the map $\mathcal{B}$ is a contraction map on $Z_{T}$ for
small $T>0$.

In fact, using that $\Vert e^{\delta t\Delta }\Vert
_{B_{p,1}^{3/p+1}\rightarrow B_{p,1}^{3/p+1}}\leq 1$ and $\mathbb{P}$ is
bounded in $B_{p,1}^{3/p+1}$ for $1<p<\infty $, we have that
\begin{equation*}
\left\Vert \int_{0}^{t}e^{\delta (t-\tau )\Delta }\mathbb{P}\Omega
e_{3}\times u(\tau )\ d\tau \right\Vert _{{B_{p,1}^{3/p+1}}}\leq C|\Omega
|\int_{0}^{t}\Vert e_{3}\times u(\tau )\Vert _{B_{p,1}^{3/p+1}}d\tau \leq
C|\Omega |T\sup_{0\leq t\leq T}\Vert u(t)\Vert _{B_{p,1}^{3/p+1}}.
\end{equation*}%
Taking the supremum over $t\in \lbrack 0,T],$ we get a constant $C>0$ such
that
\begin{equation}
\sup_{0\leq t\leq T}\left\Vert \int_{0}^{t}e^{\delta (t-\tau )\Delta }%
\mathbb{P}\Omega e_{3}\times u(\tau )\ d\tau \right\Vert _{{B_{p,1}^{3/p+1}}%
}\leq C|\Omega |T\sup_{0\leq t\leq T}\Vert u(t)\Vert _{B_{p,1}^{3/p+1}},
\label{3.4.1}
\end{equation}%
for all $u\in C([0,T];{B_{p,1}^{3/p+1}}(\mathbb{R}^{3}))$. Similarly,
\begin{equation}
\left\Vert \int_{0}^{t}e^{\delta (t-\tau )\Delta }\mathbb{P}\Omega
e_{3}\times u(\tau )\ d\tau \right\Vert _{{B_{p,1}^{3/p+2}}}\leq C|\Omega
|\int_{0}^{t}\Vert e_{3}\times u(\tau )\Vert _{B_{p,1}^{3/p+2}}d\tau \leq
C|\Omega |\Vert u\Vert _{L^{1}(0,T;B_{p,1}^{3/p+2})}  \label{aux-ineq-10}
\end{equation}%
for all $t\in \lbrack 0,T]$. An integration of (\ref{aux-ineq-10}) over $%
[0,T]$ yields the estimate
\begin{equation}
\left\Vert \int_{0}^{t}e^{\delta (t-\tau )\Delta }\mathbb{P}\Omega
e_{3}\times u(\tau )\ d\tau \right\Vert _{L^{1}(0,T;{B_{p,1}^{3/p+2}})}\leq
C|\Omega |T\Vert u\Vert _{L^{1}(0,T;{B_{p,1}^{3/p+2}})},  \label{3.4.2}
\end{equation}%
for all $u\in L^{1}(0,T;{B_{p,1}^{3/p+2}}(\mathbb{R}^{3})).$

Next we can apply (\ref{3.4.1}), (\ref{3.4.2}) and Lemma \ref{3.5} in order
to estimate

\begin{eqnarray}
\Vert \mathcal{B}(u^{\delta })-\mathcal{B}(v^{\delta })\Vert _{Z_{T}}
&=&\Vert \int_{0}^{t}e^{\delta (t-\tau )\Delta }\mathbb{P}\Omega e_{3}\times
\left( u^{\delta }(\tau )-v^{\delta }(\tau )\right) \ d\tau  \notag \\
&&+\int_{0}^{t}e^{\delta (t-\tau )\Delta }\mathbb{P}\{\left( u^{\delta
}(\tau )-v^{\delta }(\tau )\right) \cdot \nabla \}u^{\delta }(\tau )\ d\tau
\notag \\
&&+\int_{0}^{t}e^{\delta (t-\tau )\Delta }\mathbb{P}\left( v^{\delta }(\tau
)\cdot \nabla \right) \left( u^{\delta }(\tau )-v^{\delta }(\tau )\right) \
d\tau \Vert _{Z_{T}}  \notag \\
&\leq &C_{1}|\Omega |T\Vert u^{\delta }-v^{\delta }\Vert
_{Z_{T}}+C_{2}L_{\delta ,T}\left( \Vert u^{\delta }\Vert _{Z_{T}}+\Vert
v^{\delta }\Vert _{Z_{T}}\right) \Vert u^{\delta }-v^{\delta }\Vert _{Z_{T}}
\notag \\
&\leq &\left\{ C_{1}|\Omega |T+4C_{0}C_{2}L_{\delta ,T}\Vert u_{0}\Vert
_{B_{p,1}^{3/p+1}}\right\} \Vert u^{\delta }-v^{\delta }\Vert _{Z_{T}},
\label{eq3.5}
\end{eqnarray}%
for all $u^{\delta },v^{\delta }\in Z_{T}$. Moreover, using (\ref{eq3.3})
and (\ref{eq3.5}) with $v^{\delta }=0,$ we obtain

\begin{eqnarray}
\Vert \mathcal{B}(u^{\delta })\Vert _{Z_{T}} &\leq &\Vert e^{\delta t\Delta
}u_{0}\Vert _{Z_{T}}+\Vert \mathcal{B}(u^{\delta })-\mathcal{B}(0)\Vert
_{Z_{T}}  \notag \\
&\leq &C_{0}\Vert u_{0}\Vert _{B_{p,1}^{3/p+1}}+\left\{ C_{1}|\Omega
|T+4C_{0}C_{2}L_{\delta ,T}\Vert u_{0}\Vert _{B_{p,1}^{3/p+1}}\right\} \Vert
u^{\delta }\Vert _{Z_{T}}  \notag \\
&\leq &C_{0}\Vert u_{0}\Vert _{B_{p,1}^{3/p+1}}\left\{ 1+2C_{1}|\Omega
|T+8C_{0}C_{2}L_{\delta ,T}\Vert u_{0}\Vert _{B_{p,1}^{3/p+1}}\right\} ,
\label{eq3.4}
\end{eqnarray}%
for all $u^{\delta }\in Z_{T}$. Next we choose $T=T_{\delta ,\Omega
}=T(\delta ,|\Omega |,\Vert u_{0}\Vert _{B_{p,1}^{3/p+1}})>0$ such that
\begin{equation}
2C_{1}|\Omega |T_{\delta ,\Omega }+8C_{0}C_{2}\Vert u_{0}\Vert
_{B_{p,1}^{3/p+1}}\left( T_{\delta ,\Omega }+T_{\delta ,\Omega }^{\frac{1}{2}%
}\delta ^{-\frac{1}{2}}\right) <1.  \label{eq3.6}
\end{equation}%
Inserting (\ref{eq3.6}) into (\ref{eq3.4}) and (\ref{eq3.5}), we get that $%
\mathcal{B}(Z_{T_{\delta ,\Omega }})\subset Z_{T_{\delta ,\Omega }}$ and
\begin{equation*}
\Vert \mathcal{B}(u^{\delta })-\mathcal{B}(u^{\delta })\Vert _{Z_{T_{\delta
,\Omega }}}\leq \frac{1}{2}\Vert u^{\delta }-v^{\delta }\Vert _{Z_{T_{\delta
,\Omega }}},\text{ for all }u^{\delta },v^{\delta }\in Z_{T_{\delta ,\Omega
}},
\end{equation*}%
which gives the claim. By the Banach Fixed Point Theorem, there exists a
unique solution $u^{\delta }\in Z_{T_{\delta ,\Omega }}$ for (\ref{eq3.2}).

We claim that $u^{\delta }\in Z_{T_{\delta ,\Omega }}$ is a strong solution
for (\ref{Approx-scheme}) in the class (\ref{eq3.1}). By the above estimates
and using that $u^{\delta }\in C([0,T];B_{p,1}^{3/p+1}(\mathbb{R}^{3}))\cap
L^{1}(0,T;B_{p,1}^{3/p+2}(\mathbb{R}^{3}))$, it is not difficult to see that
\begin{equation*}
\mathbb{P}\Omega e_{3}\times u^{\delta }+\mathbb{P}\left( u^{\delta }\cdot
\nabla \right) u^{\delta }\in L^{1}(0,T_{\delta ,\Omega };B_{p,1}^{3/p+1}(%
\mathbb{R}^{3}))
\end{equation*}%
and $\delta \Delta v^{\delta }\in L^{1}(0,T_{\delta ,\Omega };B_{p,1}^{3/p}(%
\mathbb{R}^{3}))$ where
\begin{equation*}
v^{\delta }(t):=-\int_{0}^{t}e^{\delta (t-\tau )\Delta }\mathbb{P}\{\Omega
e_{3}\times u^{\delta }(\tau )+\left( u^{\delta }(\tau )\cdot \nabla \right)
u^{\delta }(\tau )\}\ d\tau .
\end{equation*}%
Thus, $\partial _{t}v^{\delta }\in L^{1}(0,T_{\delta ,\Omega };B_{p,1}^{3/p}(%
\mathbb{R}^{3}))$ and then $v^{\delta }\in AC([0,T_{\delta ,\Omega
}];B_{p,1}^{3/p}(\mathbb{R}^{3})).$ Moreover, $e^{\delta t\Delta }u_{0}\in
\linebreak {\ AC([0,T_{\delta ,\Omega }];B_{p,1}^{3/p}(\mathbb{R}^{3}))}$.
By standard arguments (see Kato \cite{Kato1972} and Pazy \cite{Pazy1983}),
we obtain the desired claim. For more details see \cite{Koh2015}. The
uniqueness follows from the fact that $u^{\delta }$ is the unique solution
for (\ref{eq3.2}) in the class $Z_{T_{\delta ,\Omega }}$. \fin

In what follows, we prove that there exists $T>0$ independent of $\delta \in
(0,1)$ and $\Omega \in \mathbb{R}$ such that the solution $u^{\delta }$
exists on $[0,T]$. For that, we need some \textit{a priori} uniform estimates for $%
u^{\delta }$ in the space $B_{2,1}^{5/2}.$

\begin{prop}
\label{est-seq-apro-1}Assume that $u_{0}\in B_{2,1}^{5/2}(\mathbb{R}^{3})$
and $\nabla \cdot u_{0}=0$. There exists $T=T(\Vert u_{0}\Vert
_{B_{2,1}^{5/2}})>0$ such that (\ref{Approx-scheme}) has a unique strong
solution
\begin{equation*}
u^{\delta }\in C([0,T];B_{2,1}^{5/2}(\mathbb{R}^{3})\cap
AC([0,T];B_{2,1}^{3/2}(\mathbb{R}^{3}))
\end{equation*}%
for all $0<\delta <1$ and $\Omega \in \mathbb{R}$. Furthermore, $\{u^{\delta
}\}_{\delta \in (0,1)}$ is bounded in $C([0,T];B_{2,1}^{5/2}(\mathbb{R}%
^{3})) $.
\end{prop}

\textbf{Proof.} Applying the Littlewood-Paley operator $\Delta _{j}$ to the
equation in (\ref{Approx-scheme}), taking the $L^{2}$-norm product with $%
\Delta _{j}u^{\delta }(t)$, and using $\nabla \cdot \Delta _{j}u^{\delta }=0$
and the skew-symmetric of $e_{3}\times $, we have that
\begin{equation}
\frac{1}{2}\frac{d}{dt}\Vert \Delta _{j}u^{\delta }(t)\Vert
_{L^{2}}^{2}+\delta \langle -\Delta \Delta _{j}u^{\delta }(t),\Delta
_{j}u^{\delta }(t)\rangle _{L^{2}}=-\langle \Delta _{j}(u^{\delta }(t)\cdot
\nabla )u^{\delta }(t),\Delta _{j}u^{\delta }(t)\rangle _{L^{2}}.
\label{boundedsol0}
\end{equation}%
Notice that the second term in the right hand side of \eqref{boundedsol0} is
non-negative. So, using that
\begin{equation*}
\langle (u^{\delta }(t)\cdot \nabla )\Delta _{j}u^{\delta }(t),\Delta
_{j}u^{\delta }(t)\rangle _{L^{2}}=0
\end{equation*}%
and recalling the definition of the commutator $[u^{\delta }(t)\cdot \nabla
,\Delta _{j}]$, we get
\begin{equation*}
\frac{1}{2}\frac{d}{dt}\Vert \Delta _{j}u^{\delta }(t)\Vert _{L^{2}}^{2}\leq
\langle \lbrack u^{\delta }(t)\cdot \nabla ,\Delta _{j}]u^{\delta
}(t),\Delta _{j}u^{\delta }(t)\rangle _{L^{2}}.
\end{equation*}%
By the Cauchy-Schwarz inequality, it follows that
\begin{equation*}
\frac{d}{dt}\Vert \Delta _{j}u^{\delta }(t)\Vert _{L^{2}}\leq \Vert \lbrack
u^{\delta }(t)\cdot \nabla ,\Delta _{j}]u^{\delta }(t)\Vert _{L^{2}}.
\end{equation*}%
Multiplying by $2^{5/2j}$, applying the $l^{1}(\mathbb{Z})$-norm and Lemma %
\ref{2.4}, we can estimate
\begin{equation*}
\begin{split}
\frac{d}{dt}\Vert u^{\delta }(t)\Vert _{\dot{B}_{2,1}^{5/2}}& =\sum_{j\in
\mathbb{Z}}2^{5/2j}\frac{d}{dt}\Vert \Delta _{j}u^{\delta }(t)\Vert _{L^{2}}
\\
& \leq \sum_{j\in \mathbb{Z}}2^{5/2j}\Vert \lbrack u^{\delta }(t)\cdot
\nabla ,\Delta _{j}]u^{\delta }(t)\Vert _{L^{2}} \\
& \leq C\Vert \nabla u^{\delta }(t)\Vert _{L^{\infty }}\Vert u^{\delta
}(t)\Vert _{\dot{B}_{2,1}^{5/2}}.
\end{split}%
\end{equation*}%
By Remark \ref{remark1}, it follows that
\begin{equation}
\frac{d}{dt}\Vert u^{\delta }(t)\Vert _{\dot{B}_{2,1}^{5/2}}\leq C\Vert
u^{\delta }(t)\Vert _{\dot{B}_{2,1}^{5/2}}^{2}.  \label{bounded00}
\end{equation}%
On the other hand, taking the $L^{2}$-norm product with $u^{\delta }(t)$ in (%
\ref{Approx-scheme}), we arrive at
\begin{equation*}
\frac{1}{2}\frac{d}{dt}\Vert u^{\delta }(t)\Vert _{L^{2}}^{2}+\langle
-\delta \Delta u^{\delta }(t),u^{\delta }(t)\rangle _{L^{2}}=0.
\end{equation*}%
Above, we have used the skew-symmetric of $e_{3}\times $ and $\langle
(u^{\delta }(t)\cdot \nabla )u^{\delta }(t),u^{\delta }(t)\rangle _{L^{2}}=0$
because $\nabla \cdot u^{\delta }=0$. Then,
\begin{equation}
\frac{d}{dt}\Vert u^{\delta }(t)\Vert _{L^{2}}\leq 0.  \label{bounded01}
\end{equation}%
Denote by $\Vert \cdot \Vert _{B_{2,1}^{s}}^{\ast }$ the equivalent norm $%
\Vert \cdot \Vert _{L^{2}}+\Vert \cdot \Vert _{\dot{B}_{2,1}^{s}}$ in $%
B_{2,1}^{s}$ (see (\ref{normequivalence})). By (\ref{bounded00}) and (\ref%
{bounded01}), we have that
\begin{eqnarray}
\frac{d}{dt}\Vert u^{\delta }(t)\Vert _{B_{2,1}^{5/2}}^{\ast } &=&\frac{d}{dt%
}\left( \Vert u^{\delta }(t)\Vert _{\dot{B}_{2,1}^{5/2}}+\Vert u^{\delta
}(t)\Vert _{L^{2}}\right)   \notag \\
&\leq &\frac{d}{dt}\Vert u^{\delta }(t)\Vert _{\dot{B}_{2,1}^{5/2}}  \notag
\\
&\leq &C\Vert u^{\delta }(t)\Vert _{\dot{B}_{2,1}^{5/2}}^{2}  \notag \\
&\leq &C(\Vert u^{\delta }(t)\Vert _{B_{2,1}^{5/2}}^{\ast })^{2}.
\label{aux-ineq-100}
\end{eqnarray}%
Using (\ref{aux-ineq-100}) and that $K_{1}\Vert \cdot \Vert
_{B_{2,1}^{s}}\leq \Vert \cdot \Vert _{B_{2,1}^{s}}^{\ast }\leq K_{2}\Vert
\cdot \Vert _{B_{2,1}^{s}}$ for some $K_{1},K_{2}>0$, it follows that

\begin{equation*}
\Vert u^{\delta }(t)\Vert _{B_{2,1}^{5/2}}\leq \frac{1}{K_{1}}\Vert
u^{\delta }(t)\Vert _{B_{2,1}^{5/2}}^{\ast }\leq \frac{1}{K_{1}}\frac{\Vert
u_{0}\Vert _{B_{2,1}^{5/2}}^{\ast }}{1-C\Vert u_{0}\Vert
_{B_{2,1}^{5/2}}^{\ast }t}\ \ \text{ }\leq \frac{1}{K_{1}}\frac{K_{2}\Vert
u_{0}\Vert _{B_{2,1}^{5/2}}}{1-CK_{2}\Vert u_{0}\Vert _{B_{2,1}^{5/2}}t},
\end{equation*}%
for $0\leq t<(CK_{2}\Vert u_{0}\Vert _{B_{2,1}^{5/2}})^{-1}$. Taking $T=T(\Vert u_{0}\Vert _{B_{2,1}^{5/2}})=(2CK_{2}\Vert
u_{0}\Vert _{B_{2,1}^{5/2}})^{-1}$ and $L=2K_{2}/K_{1}$, we obtain
\begin{equation}
\Vert u^{\delta }(t)\Vert _{B_{2,1}^{5/2}}\leq L\Vert u_{0}\Vert
_{B_{2,1}^{5/2}},\ \ \text{ for all }\ \ t\in \lbrack 0,T].
\label{boundedsol}
\end{equation}

Notice that $T>0$ is independent of $\delta \in (0,1)$ and $\Omega \in
\mathbb{R}$. If $T_{\delta ,\Omega }<T$, by (\ref{eq3.6}) and (\ref%
{boundedsol}) we can take $T_{\delta ,\Omega }^{\prime }=T_{\delta ,\Omega
}^{\prime }(\Vert u_{0}\Vert _{B_{2,1}^{5/2}})>0$ small enough and solve (%
\ref{Approx-scheme}) on $[T_{\delta ,\Omega },T_{\delta ,\Omega }+T_{\delta
,\Omega }^{\prime }]$ with the initial value $u^{\delta }(T_{\delta ,\Omega
})\in B_{2,1}^{5/2}(\mathbb{R}^{3})$. It follows that the solution $%
u^{\delta }$ can be extended to the interval $[0,T_{\delta ,\Omega
}+T_{\delta ,\Omega }^{\prime }]$. Invoking again the same procedure, we can
extend $u^{\delta }$ (if necessary) to $[0,T_{\delta ,\Omega }+2T_{\delta
,\Omega }^{\prime }]$, $[0,T_{\delta ,\Omega }+3T_{\delta ,\Omega }^{\prime
}]$ and so on, and obtain a solution $u^{\delta }$ for (\ref{Approx-scheme})
on $[0,T]$ satisfying (\ref{boundedsol}). \fin

\section{Proof of Theorem \protect\ref{principal-theorem}}

In this section we prove Theorem \ref{principal-theorem} through three
subsections.

\subsection{Proof of item $(i)$}

For $0<\delta _{1}<\delta _{2}<1$, we can write
\begin{equation}
\left\{
\begin{split}
& \partial _{t}(u^{\delta _{1}}-u^{\delta _{2}})-\delta _{1}\Delta
(u^{\delta _{1}}-u^{\delta _{2}})+(\delta _{2}-\delta _{1})\Delta u^{\delta
_{2}}=-\mathbb{P}\Omega e_{3}\times (u^{\delta _{1}}-u^{\delta _{2}}) \\
& -\mathbb{P}\left\{ (u^{\delta _{1}}-u^{\delta _{2}})\cdot \nabla \right\}
u^{\delta _{1}}-\mathbb{P}(u^{\delta _{2}}\cdot \nabla )(u^{\delta
_{1}}-u^{\delta _{2}}), \\
& \nabla \cdot u^{\delta _{1}}=\nabla \cdot u^{\delta _{2}}=0, \\
& (u^{\delta _{1}}-u^{\delta _{2}})(0,x)=0.
\end{split}%
\right.  \label{5.1}
\end{equation}%
We will show that there exists a limit $u\in C([0,T];B_{2,1}^{3/2}(\mathbb{R}%
^{3}))$ such that
\begin{equation}
u^{\delta }(t)\rightarrow u(t)\ \text{in}\ B_{2,1}^{3/2}\text{ uniformly for
}t\in \lbrack 0,T].\text{ }  \label{continuousclass}
\end{equation}%
We start by obtaining estimates in ${B}_{2,1}^{3/2}$ for the difference $%
u^{\delta _{1}}-u^{\delta _{2}}$ uniformly in $[0,T]$. Computing the $L^{2}$%
-inner product of (\ref{5.1}) with $u^{\delta _{1}}-u^{\delta _{2}}$, and
afterwards using the skew-symmetry of $\left( e_{3}\times \cdot \right) $, $%
\nabla \cdot (u^{\delta _{1}}-u^{\delta _{2}})=0$, Holder inequality, and
Remark \ref{remark1}, we obtain%
\begin{equation*}
\begin{split}
\frac{1}{2}\frac{d}{dt}& \Vert (u^{\delta _{1}}-u^{\delta _{2}})(t)\Vert
_{L^{2}}^{2} \\
& \leq (\delta _{2}-\delta _{1})\Vert -\Delta u^{\delta _{2}}(t)\Vert
_{L^{2}}\Vert (u^{\delta _{1}}-u^{\delta _{2}})(t)\Vert _{L^{2}}+\Vert
\nabla u^{\delta _{1}}(t)\Vert _{L^{\infty }}\Vert (u^{\delta
_{1}}-u^{\delta _{2}})(t)\Vert _{L^{2}}^{2} \\
& \leq C\delta _{2}\Vert -\Delta u^{\delta _{2}}(t)\Vert _{L^{2}}\Vert
(u^{\delta _{1}}-u^{\delta _{2}})(t)\Vert _{L^{2}}+\Vert u^{\delta
_{1}}(t)\Vert _{B_{2,1}^{5/2}}\Vert (u^{\delta _{1}}-u^{\delta
_{2}})(t)\Vert _{L^{2}}^{2}.
\end{split}%
\end{equation*}%
Integrating over $(0,t)$ and using (\ref{boundedsol}), we arrive at the
estimate
\begin{equation}
\begin{split}
\Vert (u^{\delta _{1}}-u^{\delta _{2}})(t)\Vert _{L^{2}} &\leq C\delta
_{2}\int_{0}^{t}\Vert -\Delta u^{\delta _{2}}(\tau )\Vert _{L^{2}}\ d\tau
+C\int_{0}^{t}\Vert u^{\delta _{1}}(\tau )\Vert _{B_{2,1}^{5/2}}\Vert
(u^{\delta _{1}}-u^{\delta _{2}})(\tau )\Vert _{L^{2}}\ d\tau \\
&\leq C\delta_{2}T\Vert u^{\delta _{2}}\Vert _{L^{\infty
}(0,T;B_{2,1}^{5/2})} +C\int_{0}^{t}\Vert u^{\delta_{1}}(\tau )\Vert
_{B_{2,1}^{5/2}}\Vert (u^{\delta _{1}}-u^{\delta _{2}})(\tau )\Vert
_{L^{2}}\ d\tau \\
&\leq C\delta_{2}T\Vert u_{0}\Vert_{B_{2,1}^{5/2}} +C\Vert u_{0}\Vert
_{B_{2,1}^{5/2}}\int_{0}^{t}\Vert (u^{\delta _{1}}-u^{\delta _{2}})(\tau
)\Vert _{L^{2}}\ d\tau  \label{5.3}
\end{split}%
\end{equation}

By Gronwall inequality and (\ref{5.3}), we have that there exists $C>0$ such
that
\begin{equation*}
\Vert (u^{\delta _{1}}-u^{\delta _{2}})(t)\Vert _{L^{2}}\leq C\delta
_{2}T\Vert u_{0}\Vert _{B_{2,1}^{5/2}}\exp \{C\Vert u_{0}\Vert
_{B_{2,1}^{5/2}}t\}
\end{equation*}%
and, consequently, as $\delta _{2}\rightarrow 0^{+}$ we have
\begin{equation}
\sup_{0<t<T}\Vert (u^{\delta _{1}}-u^{\delta _{2}})(t)\Vert _{L^{2}}\leq
C\delta _{2}T\Vert u_{0}\Vert _{B_{2,1}^{5/2}}\exp \{C\Vert u_{0}\Vert
_{B_{2,1}^{5/2}}T\}\rightarrow 0.  \label{aux-conv-10}
\end{equation}%
Let $0<\theta <1$ and $s_{1},s_{2},s_{3}\geq 0$ be such $s_{3}=(1-\theta
)s_{1}+\theta s_{2}$. By Gagliardo-Nirenberg type inequality in Besov spaces
(see \cite{H-M-O-W}), we can estimate
\begin{equation}
\Vert (u^{\delta _{1}}-u^{\delta _{2}})(t)\Vert _{B_{2,1}^{s_{3}}}\leq
C\Vert (u^{\delta _{1}}-u^{\delta _{2}})(t)\Vert
_{B_{2,2}^{s_{1}}}^{1-\theta }\Vert (u^{\delta _{1}}-u^{\delta
_{2}})(t)\Vert _{B_{2,1}^{s_{2}}}^{\theta }.  \label{aux-conv-11}
\end{equation}%
Considering $s_{1}=0,$ $s_{2}=5/2$ and $s_{3}=\theta s_{2}$ in (\ref%
{aux-conv-11}), and using (\ref{boundedsol}), $B_{2,2}^{0}=L^{2}$ and (\ref%
{aux-conv-10}) , we obtain
\begin{equation*}
\Vert u^{\delta _{1}}-u^{\delta _{2}}\Vert _{L^{\infty
}(0,T;B_{2,1}^{s_{3}})}\leq C\left\Vert u_{0}\right\Vert
_{B_{2,1}^{5/2}}^{\theta }\Vert u^{\delta _{1}}-u^{\delta _{2}}\Vert
_{L^{\infty }(0,T;L^{2})}^{1-\theta }\rightarrow 0,\text{ as }\delta
_{2}\rightarrow 0^{+},
\end{equation*}%
for each fixed $\theta \in (0,1).$ Hence, by completeness and uniqueness of
the limit in the distributional sense, $u^{\delta }\rightarrow u$ in $%
L^{\infty }(0,T;B_{2,1}^{\tilde{s}})$ for all $0<\tilde{s}<5/2$. In
particular, taking $\tilde{s}=3/2$ and recalling that $u^{\delta }\in
C([0,T];B_{2,1}^{3/2}(\mathbb{R}^{3})),$ we obtain (\ref{continuousclass}).

Also, in view of (\ref{boundedsol}), it follows that $(u^{\delta })_{\delta
\in (0,1)}$ is bounded in $L^{\infty }(0,T;B_{2,1}^{5/2}(\mathbb{R}^{3}))$.
Then, we can extract a subsequence $(u^{\delta^{(j)}})_{j=1}^{\infty }$ that
converges to $u$ weakly-$\star $ in $L^{\infty }(0,T;B_{2,1}^{5/2}(\mathbb{R}%
^{3}))$. Thus we have that
\begin{equation}
u\in L^{\infty }(0,T;B_{2,1}^{5/2}(\mathbb{R}^{3}))\cap
C([0,T];B_{2,1}^{3/2}(\mathbb{R}^{3}))  \label{prop-u-1}
\end{equation}%
and
\begin{equation}
\Vert u\Vert _{L^{\infty }(0,T;B_{2,1}^{5/2})}\leq \liminf_{j\rightarrow \infty}\Vert u^{\delta^{(j)}}\Vert _{L^{\infty }(0,T;B_{2,1}^{5/2})}\leq
L\Vert u_{0}\Vert _{B_{2,1}^{5/2}}.  \label{5.10}
\end{equation}%
\qquad

Next we claim that $u$ is a solution for (\ref{prin-problem-leray}). For the
nonlinear term, by using integration by parts, Lemma \ref{inequalityholder},
Remark \ref{remark1}, (\ref{boundedsol}) and (\ref{5.10}), we can estimate
\begin{equation*}
\begin{split}
\int_{0}^{t}& \Vert \mathbb{P}\nabla \cdot \left[ u^{\delta }(\tau )\otimes
u^{\delta }(\tau )-u(\tau )\otimes u(\tau )\right] \Vert _{B_{2,1}^{3/2}}\
d\tau \\
& =\int_{0}^{t}\Vert \mathbb{P}\nabla \cdot \left[ (u^{\delta }(\tau
)-u(\tau ))\otimes u^{\delta }(\tau )+u(\tau )\otimes (u^{\delta }(\tau
)-u(\tau ))\right] \Vert _{B_{2,1}^{3/2}}\ d\tau \\
& \leq C\int_{0}^{t}\left\{ \Vert u^{\delta }(\tau )\Vert
_{B_{2,1}^{5/2}}+\Vert u(\tau )\Vert _{B_{2,1}^{5/2}}\right\} \Vert
u^{\delta }(\tau )-u(\tau )\Vert _{B_{2,1}^{3/2}}\ d\tau \\
& \leq CT\Vert u_{0}\Vert _{B_{2,1}^{5/2}}\sup_{0<t<T}\Vert u^{\delta
}(t)-u(t)\Vert _{B_{2,1}^{3/2}}\rightarrow 0,\text{ as }\delta \rightarrow
0^{+}.
\end{split}%
\end{equation*}%
which implies
\begin{equation}
\int_{0}^{t}\mathbb{P}(u^{\delta }(\tau )\cdot \nabla )u^{\delta }(\tau )\
d\tau \rightarrow \int_{0}^{t}\mathbb{P}(u(\tau )\cdot \nabla )u(\tau )\
d\tau \ \ \mbox{in}\ \ L^{\infty }((0,T);B_{2,1}^{3/2}),\text{ as }\delta
\rightarrow 0^{+}\text{.}  \label{5.11}
\end{equation}%
Also, we have that
\begin{eqnarray*}
\delta \int_{0}^{t}\Vert -\Delta u^{\delta }(\tau )\Vert _{B_{2,1}^{1/2}}\
d\tau &\leq &\delta \int_{0}^{t}\Vert u^{\delta }(\tau )\Vert
_{B_{2,1}^{5/2}}\ d\tau \\
&\leq &\delta T\Vert u^{\delta }\Vert _{L^{\infty }(0,T;B_{2,1}^{5/2})}\leq
C\delta T\Vert u_{0}\Vert _{B_{2,1}^{5/2}}\rightarrow 0
\end{eqnarray*}%
and
\begin{equation*}
\int_{0}^{t}\Vert \mathbb{P}\Omega e_{3}\times (u^{\delta }(\tau )-u(\tau
))\Vert _{B_{2,1}^{3/2}}\ d\tau \leq CT|\Omega |\sup_{0<t<T}\Vert u^{\delta
}(t)-u(t)\Vert _{B_{2,1}^{3/2}}\rightarrow 0,\text{ as }\delta \rightarrow
0^{+}.
\end{equation*}%
Then
\begin{equation}
\begin{split}
& \delta \int_{0}^{t}-\Delta u^{\delta }(\tau )\ d\tau \rightarrow 0\ \ %
\mbox{in}\ L^{\infty }((0,T);B_{2,1}^{1/2}) \\
& \int_{0}^{t}\mathbb{P}\Omega e_{3}\times u^{\delta }(\tau )\ d\tau
\rightarrow \int_{0}^{t}\mathbb{P}\Omega e_{3}\times u(\tau )\ d\tau \ \ %
\mbox{in}\ \ L^{\infty }((0,T);B_{2,1}^{3/2}),\text{ as }\delta \rightarrow
0^{+}.
\end{split}
\label{5.13}
\end{equation}%
Therefore, since $u^{\delta }$ satisfies (\ref{Approx-scheme}), we obtain
from (\ref{5.11}), (\ref{5.13}) and the continuous inclusion $B_{2,1}^{3/2}\subset B_{2,1}^{1/2}$ that
\begin{equation}
u(t)-u_{0}=\int_{0}^{t}\left\{ \mathbb{P}\Omega e_{3}\times u(\tau )+\mathbb{%
P}(u(\tau )\cdot \nabla )u(\tau )\right\} \ d\tau \ \ \ \mbox{in}\ \
B_{2,1}^{1/2}(\mathbb{R}^{3}).  \label{aux-eq-euler-1}
\end{equation}%
In view of the above estimates and (\ref{prop-u-1}), we can see that both
sides of (\ref{aux-eq-euler-1}) belong to $C([0,T];B_{2,1}^{3/2}(\mathbb{R}%
^{3}).$ Thus, equality (\ref{aux-eq-euler-1}) holds in $B_{2,1}^{3/2}(%
\mathbb{R}^{3})$ and $u\in AC([0,T];B_{2,1}^{3/2}(\mathbb{R}^{3}))\cap
L^{\infty }(0,T;B_{2,1}^{5/2}(\mathbb{R}^{3}))$ is a solution for (\ref%
{prin-problem-leray}), as claimed.

The next lemma deals with the time-continuity of solutions for (\ref%
{prin-problem-leray}). In particular, for $p=2$ it implies the
time-continuity of the solution $u$ in $B_{2,1}^{5/2}(\mathbb{R}^{3})$
obtained as limit of the approximation scheme. Notice that in fact it holds
for $1\leq p<\infty .$

\begin{lem}
\label{continuous} Let $0<T<\infty $ and $1\leq p<\infty $. If $u$ is a
solution for (\ref{prin-problem-leray}) in $L^{\infty }(0,T;B_{p,1}^{3/p+1}(%
\mathbb{R}^{3}))$ with initial velocity $u_{0}\in B_{p,1}^{3/p+1}(\mathbb{R}%
^{3})$ satisfying $\nabla \cdot u_{0}=0$, then $u\in C([0,T];B_{p,1}^{3/p+1}(%
\mathbb{R}^{3}))$.
\end{lem}

\textbf{Proof.} Firstly, by Lemma \ref{inequalityholder}, we have that $%
\partial _{t}u\in L^{\infty }(0,T;B_{p,1}^{3/p}(\mathbb{R}^{3}))$. Thus

\begin{equation*}
u\in W^{1,\infty }([0,T];B_{p,1}^{3/p}(\mathbb{R}^{3}))\subset
C([0,T];B_{p,1}^{3/p}(\mathbb{R}^{3})).
\end{equation*}

For every $k\in \mathbb{N}$, we denote $w_{k}:=S_{k}u$. We are going to
prove that the sequence $\{w_{k}\}_{k\in \mathbb{N}}$ converges to $u$ in $%
L^{\infty }(0,T;B_{p,1}^{3/p}(\mathbb{R}^{3}))$. Applying the
Littlewood-Paley operator in (\ref{prin-problem-leray}), for each $j\in
\mathbb{N}$ we obtain
\begin{equation*}
\partial _{t}\Delta _{j}u+(S_{j}u\cdot \nabla )\Delta _{j}u=(S_{j}u\cdot
\nabla )\Delta _{j}u-\Delta _{j}(u\cdot \nabla )u-\Delta _{j}\nabla p-\Omega
e_{3}\times \Delta _{j}u.
\end{equation*}%
Since $\Delta _{j}u$ is absolutely continuous on $[0,T]$ with values in $%
L^{p}(\mathbb{R}^{3})$ and $\nabla \cdot S_{j-2}u=0$, we can estimate
\begin{equation*}
\begin{split}
\Vert \Delta _{j}u(t)\Vert _{L^{p}}\leq & \Vert \Delta _{j}u_{0}\Vert
_{L^{p}}+\int_{0}^{t}\Vert \Delta _{j}\nabla p\Vert _{L^{p}}\ d\tau \\
& +\int_{0}^{t}\Vert (S_{j}u\cdot \nabla )\Delta _{j}u-\Delta _{j}(u\cdot
\nabla )u\Vert _{L^{p}}\ d\tau +\int_{0}^{t}\Vert \Omega e_{3}\times \Delta
_{j}u\Vert _{L^{p}}\ d\tau .
\end{split}%
\end{equation*}%
It follows that

\begin{equation*}
\begin{split}
\Vert u(t)-w_{k}(t)\Vert _{B_{p,1}^{3/p+1}}& \leq C\sum_{j\geq
k}2^{j(3/p+1)}\Vert \Delta _{j}u(t)\Vert _{L^{p}} \\
& \leq C\bigg {(}\sum_{j\geq k}2^{j(3/p+1)}\Vert \Delta _{j}u_{0}\Vert
_{L^{p}}+\int_{0}^{t}\sum_{j\geq k}2^{j(3/p+1)}\Vert \Delta _{j}\nabla
p\Vert _{L^{p}}\ d\tau \\
& \ \ +\int_{0}^{t}\sum_{j\geq k}2^{j(3/p+1)}\Vert (S_{j}u\cdot \nabla
)\Delta _{j}u-\Delta _{j}(u\cdot \nabla )u\Vert _{L^{p}}\ d\tau \\
& \ \ +|\Omega |\int_{0}^{t}\sum_{j\geq k}2^{j(3/p+1)}\Vert \Delta
_{j}u\Vert _{L^{p}}\ d\tau \bigg {)}.
\end{split}%
\end{equation*}

The first term in the right-hand side converges to zero as $k\rightarrow
\infty $ because $u_{0}\in B_{p,1}^{3/p+1}(\mathbb{R}^{3})$. By Lemma \ref%
{inequalityholder}, Lemma \ref{bilinearest2} and the fact that $u(t)\in
B_{p,1}^{3/p+1}(\mathbb{R}^{3})$, we have that the second and third terms in
the right-hand side also converge to zero as $k\rightarrow \infty $.
Therefore, the sequence $\{w_{k}\}_{k\in \mathbb{N}}$ converges to $u$ in $%
L^{\infty }(0,T;B_{p,1}^{3/p+1}(\mathbb{R}^{3}))$. Moreover, we get
\begin{equation}
\begin{split}
\Vert w_{k}(s)-w_{k}(t)\Vert _{B_{p,1}^{3/p+1}}& =\Vert
S_{k}(u(s)-u(t))\Vert _{B_{p,1}^{3/p+1}} \\
& \leq C\sum_{j=-1}^{k+1}2^{j(3/p+1)}\Vert \Delta _{j}(u(s)-u(t))\Vert
_{L^{p}} \\
& \leq C2^{k+1}\Vert u(s)-u(t)\Vert _{B_{p,1}^{3/p}}.
\end{split}
\label{aux-reg-time-1}
\end{equation}

Estimate (\ref{aux-reg-time-1}) and the fact that $u\in
C([0,T];B_{p,1}^{3/p}(\mathbb{R}^{3}))$ imply that each $w_{k}\in
C([0,T];B_{p,1}^{3/p+1}(\mathbb{R}^{3}))$. Therefore, the limit $u$ also
belongs to $C([0,T];B_{p,1}^{3/p+1}(\mathbb{R}^{3}))$. \fin

Now, taking $p=2$ in Lemma \ref{continuous}, since $u\in L^{\infty
}(0,T;B_{2,1}^{5/2}(\mathbb{R}^{3}))$ and $u_{0}\in B_{2,1}^{5/2}(\mathbb{R}%
^{3})$ we have that $u\in C([0,T];B_{2,1}^{5/2}(\mathbb{R}^{3}))$, and then $%
u$ satisfies
\begin{equation}
\partial _{t}u=-\mathbb{P}\Omega e_{3}\times u-\mathbb{P}(u\cdot \nabla
)u\in C([0,T];B_{2,1}^{3/2}(\mathbb{R}^{3})).  \label{deriv-time}
\end{equation}%
This shows that $u\in C^{1}([0,T];B_{2,1}^{3/2}(\mathbb{R}^{3}))$, and
therefore $u$ is a strong solution for (\ref{prin-problem-leray}) in the
class
\begin{equation}
C([0,T];B_{2,1}^{5/2}(\mathbb{R}^{3}))\cap C^{1}([0,T];B_{2,1}^{3/2}(\mathbb{%
R}^{3})).  \label{class}
\end{equation}

\textbf{Uniqueness. }Let $u$ and $v$ be strong solutions for (\ref%
{prin-problem-leray}) in the class (\ref{class}) with the same initial data $%
u_{0}(x)$. Subtracting the corresponding equations satisfied by $u$ and $v$,
we get
\begin{equation}
\left\{
\begin{split}
& \partial _{t}(u-v)+\mathbb{P}\Omega e_{3}\times (u-v)+\mathbb{P}%
\{(u-v)\cdot \nabla \}u+\mathbb{P}(v\cdot \nabla )(u-v)=0, \\
& \nabla \cdot u=\nabla \cdot v=0, \\
& (u-v)(0,x)=0.
\end{split}%
\right.  \label{5.15}
\end{equation}%
Computing the $L^{2}$-inner product of (\ref{5.15}) with $u-v$, we obtain
\begin{equation*}
\begin{split}
\frac{1}{2}\frac{d}{dt}\Vert (u-v)(t)\Vert _{L^{2}}^{2}& =-\langle (u-v)(t),%
\mathbb{P}\{(u-v)(t)\cdot \nabla \}u(t)\rangle _{L^{2}} \\
& \leq \Vert \nabla u(t)\Vert _{L^{\infty }}\Vert (u-v)(t)\Vert _{L^{2}}^{2}
\\
& \leq C\Vert u(t)\Vert _{B_{2,1}^{5/2}}\Vert (u-v)(t)\Vert _{L^{2}}^{2},
\end{split}%
\end{equation*}%
and then

\begin{eqnarray}
\Vert (u-v)(t)\Vert _{L^{2}} &\leq &C\int_{0}^{t}\Vert u(\tau )\Vert
_{B_{2,1}^{5/2}}\Vert (u-v)(\tau )\Vert _{L^{2}}\ d\tau  \notag \\
&\leq &C\Vert u\Vert _{L^{\infty }(0,T;B_{2,1}^{5/2})}\int_{0}^{t}\Vert
(u-v)(\tau )\Vert _{L^{2}}\ d\tau .  \label{5.16}
\end{eqnarray}%
Since $\Vert u\Vert _{L^{\infty }(0,T;B_{2,1}^{5/2})}<\infty ,$ we can use
Gronwall inequality to obtain $\Vert u(t)-v(t)\Vert _{L^{2}}=0$ for all $%
t\in \lbrack 0,T]$, and then $u\equiv v.$\fin

\subsection{Blow-up criterion}

In this part, we prove a blow-up criterion of BKM type (see \cite{BKM}). We
will use it to prove item (ii) of Theorem \ref{principal-theorem}.

\begin{proposition}
\label{blowup} Let $u_{0}\in B_{2,1}^{5/2}(\mathbb{R}^{3})$ with $\nabla
\cdot u_{0}=0$. Assume that
\begin{equation}
u\in C([0,T);B_{2,1}^{5/2}(\mathbb{R}^{3}))\cap C^{1}([0,T);B_{2,1}^{3/2}(%
\mathbb{R}^{3}))  \label{aux-blow-up-1}
\end{equation}%
is a solution for (\ref{prin-problem-leray}). For some $T^{\prime }>T,$ $u$
can be extended to $[0,T^{\prime })$ with $u\in C([0,T^{\prime
});B_{2,1}^{5/2}(\mathbb{R}^{3}))\cap \linebreak {C^{1}([0,T^{\prime
});B_{2,1}^{3/2}(\mathbb{R}^{3}))}$ provided that $\int_{0}^{T}\Vert \nabla
u(t)\Vert _{L^{\infty }}\ dt<\infty $.
\end{proposition}

\textbf{Proof.} Item (i) of Theorem \ref{principal-theorem} assures that the
existence-time $T>0$ depends only on the initial data norm $\Vert u_{0}\Vert
_{B_{2,1}^{5/2}}$. Computing the $L^{2}$-inner product of (\ref%
{prin-problem-leray}) with $u$, using the symmetry of $e_{3}\times u$ and $%
\nabla \cdot u=0$, one can deduce%
\begin{equation}
\Vert u(t)\Vert _{L^{2}}=\Vert u_{0}\Vert _{L^{2}}\ \ \mbox{for all}\ \ t\in
\lbrack 0,T).  \label{6.1}
\end{equation}%
Moreover, we can apply the operator $\Delta _{j}$ in (\ref%
{prin-problem-leray}), multiply the result by $\Delta _{j}u$ and after use $%
\langle (u\cdot \nabla )\Delta _{j}u,\Delta _{j}u\rangle _{L^{2}}=0$ to get
the identity
\begin{equation}
\frac{1}{2}\frac{d}{dt}\Vert \Delta _{j}u(t)\Vert _{L^{2}}^{2}=-\langle
\Delta _{j}(u(t)\cdot \nabla )u(t),\Delta _{j}u(t)\rangle _{L^{2}}=\langle
\lbrack u(t)\cdot \nabla ,\Delta _{j}]u(t),\Delta _{j}u(t)\rangle _{L^{2}}.
\label{aux-30}
\end{equation}%
Using the Schwartz inequality and integrating (\ref{aux-30}) over $(0,t)$,
we obtain
\begin{equation}
\Vert \Delta _{j}u(t)\Vert _{L^{2}}\leq \Vert \Delta _{j}u_{0}\Vert
_{L^{2}}+\int_{0}^{t}\Vert \lbrack u(\tau )\cdot \nabla ,\Delta _{j}]u(\tau
)\Vert _{L^{2}}\ d\tau .  \label{aux-31}
\end{equation}%
Now we multiply (\ref{aux-31}) by $2^{(5/2)j}$ and afterwards take the $%
l^{1}(\mathbb{Z})$-norm to deduce
\begin{equation*}
\Vert u(t)\Vert _{{\dot{B}}_{2,1}^{5/2}}\leq \Vert u_{0}\Vert _{{\dot{B}}%
_{2,1}^{5/2}}+\int_{0}^{t}\sum_{j\in \mathbb{Z}}2^{(5/2)j}\Vert \lbrack
u(\tau )\cdot \nabla ,\Delta _{j}]u(\tau )\Vert _{L^{2}}\ d\tau .
\end{equation*}%
By Lemma \ref{2.4} (i), there exists $C>0$ such that
\begin{equation}
\Vert u(t)\Vert _{{\dot{B}}_{2,1}^{5/2}}\leq \Vert u_{0}\Vert _{{\dot{B}}%
_{2,1}^{5/2}}+C\int_{0}^{t}\Vert \nabla u(\tau )\Vert _{L^{\infty }}\Vert
u(\tau )\Vert _{B_{2,1}^{5/2}}\ d\tau .  \label{aux-32}
\end{equation}%
Putting together (\ref{6.1}) and (\ref{aux-32}), we have that
\begin{equation*}
\Vert u(t)\Vert _{B_{2,1}^{5/2}}\leq C_{3}\Vert u_{0}\Vert
_{B_{2,1}^{5/2}}+C_{4}\int_{0}^{t}\Vert \nabla u(\tau )\Vert _{L^{\infty
}}\Vert u(\tau )\Vert _{B_{2,1}^{5/2}}\ d\tau ,
\end{equation*}%
where $C_{3}$ and $C_{4}$ are positive constants. By Gronwall inequality, we
get
\begin{equation}
\Vert u(t)\Vert _{B_{2,1}^{5/2}}\leq C_{3}\Vert u_{0}\Vert _{B_{2,1}^{5/2}}%
\mbox{exp}\left\{ C_{4}\int_{0}^{t}\Vert \nabla u(\tau )\Vert _{L^{\infty
}}\ d\tau \right\} ,\text{ for all }t\in \lbrack 0,T).  \label{6.2}
\end{equation}%
Therefore, by standard arguments, if $\int_{0}^{T}\Vert \nabla u(t)\Vert
_{L^{\infty }}\ dt<\infty $ then $u$ can be continued to $[0,T]$ and so to $%
[0,T^{\prime })$ for some $T^{\prime }>T$ (by item (i) of Theorem \ref%
{principal-theorem}). \fin

The contrapositive assertion of Proposition \ref{blowup} gives the following
remark.

\begin{obs}
\label{rem-blow}Let $u_{0}\in B_{2,1}^{5/2}(\mathbb{R}^{3})$ with $\nabla
\cdot u_{0}=0$. Assume that $u$ is a solution for (\ref{prin-problem-leray})
in the class (\ref{aux-blow-up-1}). If $T=T^{\ast }<\infty$ is the maximal
existence-time, then
\begin{equation*}
\int_{0}^{T^{\ast }}\Vert \nabla u(t)\Vert _{L^{\infty }}\ dt=\infty .
\end{equation*}
\end{obs}

\subsection{Proof of item $(ii)$}

Let $u_{0}\in B_{2,1}^{7/2}(\mathbb{R}^{3})$ with $\nabla \cdot u_{0}=0$ and
let $u\in C([0,T_{\ast });B_{2,1}^{7/2}(\mathbb{R}^{3}))\cap
C^{1}([0,T_{\ast });B_{2,1}^{5/2}(\mathbb{R}^{3}))$ be the solution of (\ref%
{prin-problem-leray}) with maximal existence-time $T_{\ast }>0$. Applying
the projection operators $P_{\pm }$ in (\ref{prin-problem-leray}), we get

\begin{equation*}
\partial _{t}P_{\pm }u\mp i\Omega \frac{D_{3}}{|D|}P_{\pm }u+P_{\pm }(u\cdot
\nabla )u=0\text{ \ with }P_{\pm }u(0,x)=P_{\pm }u_{0}.
\end{equation*}%
Denoting $A_{\pm }:=\pm i\Omega \frac{D_{3}}{|D|}$ and using Duhamel principle, we have that
\begin{equation}
P_{\pm }u(t)=e^{\pm i\Omega t\frac{D_{3}}{|D|}}P_{\pm
}u_{0}-\int_{0}^{t}e^{\pm i\Omega (t-\tau )\frac{D_{3}}{|D|}}P_{\pm }(u(\tau
)\cdot \nabla )u(\tau )\ d\tau .  \label{integral}
\end{equation}

Before proceeding, we recall the Strichartz estimates of \cite{Koh2015}
which states that if $2\leq r,\theta \leq \infty $ with $(r,\theta )\neq
(2,\infty )$ and $\frac{1}{r}+\frac{1}{\theta }\leq \frac{1}{2}$ then
\begin{equation}
\Vert e^{\pm it\frac{D_{3}}{|D|}}f\Vert _{L^{r}(0,\infty ;L^{\theta })}\leq
C\Vert f\Vert _{L^{2}}.  \label{Strichartz}
\end{equation}%
Let $2<r<\infty .$ A scaling argument in (\ref{Strichartz}) leads us to
\begin{equation}
\Vert \Delta _{j}e^{\pm i\Omega t\frac{D_{3}}{|D|}}f\Vert _{L^{r}(0,\infty
;L^{\infty })}\leq C2^{\frac{3}{2}j}|\Omega |^{-\frac{1}{r}}\Vert \Delta
_{j}f\Vert _{L^{2}},  \label{8.3}
\end{equation}%
for all $j\in \mathbb{Z}$ and $\Omega \in \mathbb{R}\setminus \{0\}$, where $%
C=C(r)$ is a constant.

In what follows, we derive an estimate in $B_{\infty ,1}^{1}$ for the
solution $u$. Using $u=P_{+}u+P_{-}u$ (see Lemma \ref{projection}), we only
need to show the estimate for $P_{+}u$ and $P_{-}u$. First notice that
\begin{equation*}
\begin{split}
\Vert e^{\pm i\Omega t\frac{D_{3}}{|D|}}P_{\pm }u_{0}\Vert _{L^{r}(0,\infty ;%
{\dot{B}}_{\infty ,1}^{1})}& =\left\Vert \sum_{j\in \mathbb{Z}}2^{j}\Vert
\Delta _{j}e^{\pm i\Omega t\frac{D_{3}}{|D|}}P_{\pm }u_{0}\Vert _{L^{\infty
}}\right\Vert _{L_{t}^{r}(0,\infty )} \\
& \leq \sum_{j\in \mathbb{Z}}2^{j}\Vert \Delta _{j}e^{\pm i\Omega t\frac{%
D_{3}}{|D|}}P_{\pm }u_{0}\Vert _{L_{t}^{r}(0,\infty ;L^{\infty })} \\
& \leq C|\Omega |^{-\frac{1}{r}}\sum_{j\in \mathbb{Z}}2^{j}(2^{j})^{3/2}%
\Vert \Delta _{j}P_{\pm }u_{0}\Vert _{L^{2}} \\
& =C|\Omega |^{-\frac{1}{r}}\Vert P_{\pm }u_{0}\Vert _{{\dot{B}}%
_{2,1}^{5/2}}.
\end{split}%
\end{equation*}%
Moreover, by (\ref{Strichartz}), we have that
\begin{equation*}
\Vert e^{\pm i\Omega t\frac{D_{3}}{|D|}}P_{\pm }u_{0}\Vert _{L^{r}(0,\infty
;L^{\infty })}\leq C|\Omega |^{-\frac{1}{r}}\Vert P_{\pm }u_{0}\Vert
_{L^{2}}.
\end{equation*}%
For $2<r<\infty $ and $\Omega \in \mathbb{R}\setminus \{0\}$, the last two
inequalities yield
\begin{equation}
\Vert e^{\pm i\Omega t\frac{D_{3}}{|D|}}P_{\pm }u_{0}\Vert _{L^{r}(0,\infty ;%
{B}_{\infty ,1}^{1})}\leq C|\Omega |^{-\frac{1}{r}}\Vert P_{\pm }u_{0}\Vert
_{{B}_{2,1}^{5/2}}.  \label{5.23}
\end{equation}

For the nonlinear term, using similar arguments we obtain
\begin{equation*}
\begin{split}
\left\Vert \int_{0}^{t}e^{\pm i\Omega (t-\tau )\frac{D_{3}}{|D|}}P_{\pm
}(u(\tau )\cdot \nabla )u(\tau )\ d\tau \right\Vert _{L^{r}(0,T;L^{\infty
})}& \leq C|\Omega |^{-\frac{1}{r}}\int_{0}^{T}\Vert P_{\pm }(u(\tau )\cdot
\nabla )u(\tau )\Vert _{L^{2}}\ d\tau , \\
\left\Vert \int_{0}^{t}e^{\pm i\Omega (t-\tau )\frac{D_{3}}{|D|}}P_{\pm
}(u(\tau )\cdot \nabla )u(\tau )\ d\tau \right\Vert _{L^{r}(0,T;{\dot{B}}%
_{\infty ,1}^{1})}& \leq C|\Omega |^{-\frac{1}{r}}\int_{0}^{T}\Vert P_{\pm
}(u(\tau )\cdot \nabla )u(\tau )\Vert _{{\dot{B}}_{2,1}^{\frac{5}{2}}}\
d\tau .
\end{split}%
\end{equation*}%
Therefore
\begin{equation}
\left\Vert \int_{0}^{t}e^{\pm i\Omega (t-\tau )\frac{D_{3}}{|D|}}(u(\tau
)\cdot \nabla )u(\tau )\ d\tau \right\Vert _{L^{r}(0,T;B_{\infty
,1}^{1})}\leq C|\Omega |^{-\frac{1}{r}}\int_{0}^{T}\Vert (u(\tau )\cdot
\nabla )u(\tau )\Vert _{B_{2,1}^{\frac{5}{2}}}\ d\tau .  \label{5.24}
\end{equation}%
Estimates (\ref{5.23}) and (\ref{5.24}) imply that
\begin{equation}
\Vert u\Vert _{L^{r}(0,T;B_{\infty ,1}^{1})}\leq C|\Omega |^{-\frac{1}{r}%
}\left( \Vert u_{0}\Vert _{B_{2,1}^{\frac{5}{2}}}+\int_{0}^{T}\Vert (u(\tau
)\cdot \nabla )u(\tau )\Vert _{B_{2,1}^{\frac{5}{2}}}\ d\tau \right) ,
\label{criticalinequality}
\end{equation}%
for all $0<T<T_{\ast }$. Next, we define
\begin{equation*}
U(t):=\int_{0}^{t}\Vert \nabla u(\tau )\Vert _{L^{\infty }}\ d\tau ,\text{ \
for }0\leq t\leq T_{\ast }.
\end{equation*}%
Using the embedding $B_{\infty ,1}^{1}(\mathbb{R}^{3})\hookrightarrow
W^{1,\infty }(\mathbb{R}^{3})$, (\ref{6.2}) and estimate (\ref%
{criticalinequality}), we obtain

\begin{equation*}
\begin{split}
U(t)& \leq \int_{0}^{t}\Vert u(\tau )\Vert _{B_{\infty ,1}^{1}}\ d\tau
\\
& \leq Ct^{1-\frac{1}{r}}\Vert u\Vert _{L^{r}(0,t;B_{\infty ,1}^{1})} \\
& \leq Ct^{1-\frac{1}{r}}|\Omega |^{-\frac{1}{r}}\left( \Vert u_{0}\Vert
_{B_{2,1}^{\frac{5}{2}}}+\int_{0}^{t}\Vert (u(\tau )\cdot \nabla )u(\tau
)\Vert _{B_{2,1}^{\frac{5}{2}}}\ d\tau \right) \\
& \leq Ct^{1-\frac{1}{r}}|\Omega |^{-\frac{1}{r}}\left( \Vert u_{0}\Vert
_{B_{2,1}^{\frac{7}{2}}}+\int_{0}^{t}\Vert u(\tau )\Vert _{B_{2,1}^{\frac{7}{%
2}}}^{2}\ d\tau \right) \\
& \leq Ct^{1-\frac{1}{r}}|\Omega |^{-\frac{1}{r}}\left( \Vert u_{0}\Vert
_{B_{2,1}^{\frac{7}{2}}}+\Vert u_{0}\Vert _{B_{2,1}^{\frac{7}{2}%
}}^{2}\int_{0}^{t}\exp (CU(\tau ))d\tau \right) .
\end{split}%
\end{equation*}%
Then, there exist positive constants $C_{5}$ and $C_{6}$ (independent of $%
\Omega $) such that
\begin{equation}
U(t)\leq C_{5}t^{1-\frac{1}{r}}|\Omega |^{-\frac{1}{r}}\Vert u_{0}\Vert
_{B_{2,1}^{\frac{7}{2}}}\left( 1+\Vert u_{0}\Vert _{B_{2,1}^{7/2}}t\exp
(C_{6}U(t))\right) ,\text{ }\forall t\in (0,T_{\ast }).  \label{8.10}
\end{equation}%
For $0<T<\infty $, we consider
\begin{equation*}
H_{T}=\{t\in \lbrack 0,T]\cap \lbrack 0,T_{\ast })\mid U(t)\leq C_{5}T^{1-%
\frac{1}{r}}\Vert u_{0}\Vert _{B_{2,1}^{7/2}}\},\ \ \widehat{T}_{\ast }=\sup
H_{T}.
\end{equation*}%
We will show that $\widehat{T}_{\ast }=\min \{T,T_{\ast }\}$. For that,
suppose that $\widehat{T}_{\ast }<\min \{T,T_{\ast }\}$ by contradiction.
Then there exists $\widehat{T}$ such that $\widehat{T}_{\ast }<\widehat{T}%
<\min \{T,T_{\ast }\}$. In view of $u\in C([0,\widehat{T}];B_{2,1}^{7/2}(%
\mathbb{R}^{3}))$, we have that $U(t)$ is uniformly continuous on $[0,%
\widehat{T}]$ and
\begin{equation}
U(\widehat{T}_{\ast })\leq C_{5}T^{1-\frac{1}{r}}\Vert u_{0}\Vert
_{B_{2,1}^{7/2}}.  \label{8.11}
\end{equation}%
Taking a sufficiently large $\Omega \in \mathbb{R}\setminus \{0\}$ in such a
way that
\begin{equation}
|\Omega |^{\frac{1}{r}}\geq 2\left( 1+\Vert u_{0}\Vert _{B_{2,1}^{7/2}}T\exp
(C_{5}C_{6}T^{1-\frac{1}{r}}\Vert u_{0}\Vert _{B_{2,1}^{7/2}})\right) ,
\label{8.12}
\end{equation}%
and using (\ref{8.10}), (\ref{8.11}) and (\ref{8.12}), it follows that

\begin{equation*}
\begin{split}
U(\widehat{T}_{\ast })& \leq C_{5}(\widehat{T}_{\ast })^{1-\frac{1}{r}%
}|\Omega |^{-\frac{1}{r}}\Vert u_{0}\Vert _{B_{2,1}^{7/2}}\left( 1+\Vert
u_{0}\Vert _{B_{2,q}^{7/2}}\widehat{T}_{\ast }\exp (C_{6}U(\widehat{T}_{\ast
}))\right) \\
& \leq C_{5}T^{1-\frac{1}{r}}\Vert u_{0}\Vert _{B_{2,1}^{7/2}}|\Omega |^{-%
\frac{1}{r}}\left( 1+\Vert u_{0}\Vert _{B_{2,1}^{7/2}}T\exp (C_{5}C_{6}T^{1-%
\frac{1}{r}}\Vert u_{0}\Vert _{B_{2,1}^{7/2}})\right) \\
& \leq \frac{1}{2}C_{5}T^{1-\frac{1}{r}}\Vert u_{0}\Vert _{B_{2,1}^{7/2}}.
\end{split}%
\end{equation*}%
Thus, there exists $L$ such that $\widehat{T}_{\ast }<L<\widehat{T}$ and $%
U(L)\leq C_{5}T^{1-\frac{1}{r}}\Vert u_{0}\Vert _{B_{2,1}^{7/2}}$,
contradicting the definition of $\widehat{T}_{\ast }$. Therefore, if (\ref%
{8.12}) holds true we have that $\widehat{T}_{\ast }=\min \{T,T_{\ast }\}.$
If $T_{\ast }<T$, it follows that $T_{\ast }=\widehat{T}_{\ast }=\sup H_{T}$
and then
\begin{equation*}
U(t)=\int_{0}^{t}\Vert \nabla u(\tau )\Vert _{L^{\infty }}\ d\tau \leq
C_{5}T^{1-\frac{1}{r}}\Vert u_{0}\Vert _{B_{2,1}^{7/2}}<\infty ,
\end{equation*}%
for all $0\leq t<T_{\ast }$, and so $U(T_{\ast})<\infty$. In view of the blow-up criterion (see Remark %
\ref{rem-blow}), we are done. \fin

\end{document}